\documentclass{compositio}

\input xypic
\usepackage{amssymb,amsmath}
\usepackage{amsfonts}
\usepackage{latexsym}
\usepackage{amscd}
\usepackage{epsfig}
\usepackage{euler, latexsym, epsfig, epic}

\newtheorem{corollary}[equation]{Corollary} 
\newtheorem{lemma}[equation]{Lemma} 
\newtheorem{proposition}[equation]{Proposition} 
\newtheorem{theorem}[equation]{Theorem} 

\theoremstyle{definition} 
 
\newtheorem{definition}[equation]{Definition}

\numberwithin{equation}{section}

\begin{document}
\def\cox{\operatorname{cox}}
\def\camb{\mathfrak{C}}
\def\soc{\operatorname{soc}}
\def\rep{\operatorname{rep}}
\def\dim{\operatorname{dim}}
\def\ch{\operatorname{ch}}
\def\td{\operatorname{td}}
\def\mod{\operatorname{mod}}
\def\grmod{\operatorname{grmod}}
\def\im{\operatorname{im}}
\def\A{\mathcal{A}}
\def\add{\operatorname{add}}
\def\I{\mathcal{I}}
\def\inf{\operatorname{inf}}
\def\discrep{\operatorname{discrep}}
\def\subs{\operatorname{subs}}
\def\w{\omega}
\def\O{\mathcal{O}}
\def\pc{\mathcal{T}}
\def\u{\mu}
\def\Gen{\operatorname{Gen}}
\def\Cogen{\operatorname{Cogen}}
\def\sub{\subset}
\def\basic{\operatorname{bsc}}
\def\T{\Theta}
\def\m{\mathfrak{m}}
\def\F{\mathcal{F}}
\def\E{\mathcal{E}}
\def\ab{\operatorname{\mathfrak{a}}}
\def\faces{\mathcal{F}}
\def\P{\mathcal{P}}
\def\N{\mathcal{N}}
\def\R{\mathbb{R}}
\def\Z{\mathbb{Z}}
\def\Q{\mathbb{Q}}
\def\dimv{\operatorname{\underline{\dim}}}
\def\ww{\wedge}
\def\a{\alpha}
\def\b{\beta}
\def\c{\gamma}
\def\d{\delta}
\def\D{\mathcal{D}}
\def\e{\varepsilon}
\def\l{\lambda}
\def\O{\mathcal{O}}
\def\t{\tau}
\def\s{\sigma}
\def\di{\partial}
\def\z{\zeta}
\def\la{\langle} 
\def\ra{\rangle}
\def\rtar{\rightarrow}
\def\faces{\operatorname{faces}}
\def\can{\operatorname{can}}
\def\Ind{\operatorname{Ind}}
\def\rank{\operatorname{rank}}
\def\pic{\operatorname{pic}}
\def\symPic{\operatorname{symPic}}
\def\H{\operatorname{H}}
\def\coh{\operatorname{coh}}
\def\udim{\operatorname{udim}}
\def\order{\operatorname{ord}}
\def\height{\operatorname{ht}}
\def\det{\operatorname{det}}
\def\GL{\operatorname{GL}}
\def\PSL{\operatorname{PSL}}
\def\rad{\operatorname{rad}}
\def\biDiv{\operatorname{biDiv}}
\def\biPic{\operatorname{biPic}}
\def\locPic{\operatorname{locPic}}
\def\rtr{\operatorname{rtr}}
\def\Stan{\operatorname{Stan}}
\def\Trg{\operatorname{Trg}}
\def\B{\mathcal{B}}
\def\Sc{\mathcal{S}}
\def\Cl{\operatorname{Cl}}
\def\H{\mathcal{H}}
\def\T{\mathcal{T}}
\def\ann{\operatorname{ann}}
\def\M{\mathcal{M}}
\def\C{\mathcal{C}}
\def\CL{\mathcal{CC}}
\def\cok{\operatorname{cok}}
\def\Mod{\operatorname{Mod}}
\def\E{\mathcal{E}}
\def\ker{\operatorname{ker}}
\def\Hom{\operatorname{Hom}}
\def\shom{\mathcal{H}om}
\def\Set{\operatorname{Set}}
\def\Grp{\operatorname{Grp}}
\def\Ext{\operatorname{Ext}}
\def\Tor{\operatorname{Tor}}
\def\tors{\operatorname{tors}}
\def\End{\operatorname{End}}
\def\Pic{\operatorname{Pic}}
\def\Spec{\operatorname{Spec}}
\def\Proj{\operatorname{Proj}}
\def\Der{\operatorname{Der}}
\def\Sing{\operatorname{Sing}}
\def\OutDer{\operatorname{OutDer}}
\def\Obs{\operatorname{Obs}}
\def\AlgExt{\operatorname{AlgExt}}
\def\Exal{\operatorname{AlgExt}}
\def\Def{\operatorname{AlgExt}}
\def\kod{\operatorname{kod}}
\def\Open{\operatorname{Open}}
\def\supp{\operatorname{supp}}
\def\Br{\operatorname{Br}}
\def\EmbDef{\operatorname{EmbDef}}
\def\Aut{\operatorname{Aut}}
\def\Int{\operatorname{Int}}
\def\depth{\operatorname{depth}}
\def\lt{\operatorname{lt}}
\def\cal{\mathcal}
\def\L{\mathcal{L}}
\def\Sym{\operatorname{Sym}}
\def\PGL{\operatorname{PGL}}
\def\rk{\operatorname{rk}}
\def\Obs{\operatorname{Obs}}
\def\Out{\operatorname{Out}}
\def\ideal{\operatorname{I}}
\def\cir{\cdot}
\def\gkdim{\operatorname{gkdim}}
\def\GL{\operatorname{GL}}
\def\SL{\operatorname{SL}}
\def\Irr{\operatorname{Irr}}
\def\qcoh{\operatorname{qcoh}}
\def\id{\operatorname{id}}
\def\projd{\operatorname{pd}}
\def\G{\mathbb{G}}
\def\Target{\operatorname{Target}}
\def\Morph{\operatorname{Morph}}
\def\Ob{\operatorname{Ob}}
\def\cEnd{\operatorname{\mathcal{E}nd\,}}

\newcommand{\pt}{\mathcal{T}}
\newcommand{\ps}{\mathcal{S}}
\newcommand{\pv}{\mathcal{V}}
\newcommand{\pa}{\mathcal{A}}
\newcommand{\pb}{\mathcal{B}}

\newcommand{\pd}{\mathcal{D}}
\newcommand{\pe}{\mathcal{E}}
\newcommand{\prf}{\mathcal{F}}
\newcommand{\uperp}{{}^\perp}
\newcommand{\semi}{\semidirect}
\newcommand{\ind}{\operatorname{ind}}
\newcommand{\NC}{\operatorname{NC}}
\newcommand{\nc}{\operatorname{NC}}
\newcommand{\W}{\mathcal{W}}
\newcommand{\br}{\operatorname{Br}}
\newcommand{\Wd}{\mathcal{W}}

\newcommand{\sort}{\operatorname{sort}}
\newcommand{\ncr}{\operatorname{nc}}
\newcommand{\cl}{\operatorname{cl}}

\long\def\symbolfootnote[#1]#2{\begingroup%
\def\thefootnote{\fnsymbol{footnote}}\footnote[#1]{#2}\endgroup}

\title{Noncrossing partitions and representations of quivers}
\author{Colin Ingalls}
\email{cingalls@unb.ca}
\address{Department of Mathematics and Statistics, 
University of New Brunswick,
Fredericton, New Brunswick E3B 5A3,
Canada.}
\author{Hugh Thomas}
\email{hthomas@unb.ca}
\address{Department of Mathematics and Statistics, 
University of New Brunswick,
Fredericton, New Brunswick E3B 5A3,
Canada.}

\thanks{Both authors were supported by NSERC Discovery Grants.}

\begin{abstract}
We situate the noncrossing partitions associated to a finite Coxeter 
group within the
context of the representation theory of quivers.  We describe Reading's 
bijection between noncrossing partitions and clusters in this context,
and show that it
extends to the extended Dynkin case.  Our setup also yields
a new proof that the noncrossing partitions associated to a finite Coxeter
group form a lattice.  

We also prove some new results within the theory of quiver representations.
We show that the finitely generated, exact abelian,
and extension-closed subcategories of the representations of a quiver 
$Q$ without oriented cycles are in natural bijection with the cluster tilting
objects in the associated cluster category.  We also show these subcategories
are exactly the finitely generated categories that can be obtained as 
the semistable objects with respect to some stability condition.
\end{abstract}

\maketitle

\section{Introduction}

A partially ordered set called the
{\it noncrossing partitions} 
of $\{1,\dots,n\}$ was introduced by Kreweras \cite{Kr}
in 1972.  It was later recognized that these noncrossing partitions 
should be considered to be connected to the Coxeter group of type 
$A_{n-1}$ (that is, the symmetric group $S_n$).  
In 1997, a version of noncrossing partitions associated to type
$B_n$ was introduced by Reiner \cite{Rei}.  
The definition of noncrossing partitions for an arbitrary Coxeter group
was apparently a part of folklore before it was
written down shortly thereafter \cite{BW,Be}.  

Subsequently, {\it cluster algebras} were developed by
Fomin and Zelevinsky \cite{FZ1}.  
A cluster algebra has a set of distinguished
generators grouped into overlapping sets called {\it clusters}.
It was observed \cite{FZ} that the number of clusters for the cluster algebra 
associated to a certain
orientation of a Dynkin diagram was the same as the number of noncrossing
partitions, the generalized Catalan number.
  The reason for this was not at all obvious, though somewhat
intricate bijections have since been found \cite{Re2, ABMW}.  

The representation theory of hereditary algebras has proved an extremely 
fruitful perspective on cluster algebras from \cite{MRZ,BMRRT} 
to the more recent 
\cite{CK1,CK2}.  In this context, clusters appear as the cluster tilting objects in
the cluster category.  We will adopt this perspective on clusters throughout
this paper. 

Our goal in this paper is 
to apply the representation theory of hereditary algebras to 
account for and generalize 
two properties of the noncrossing partitions in finite type:
\begin{enumerate}
\item\label{enum} The already-mentioned fact that noncrossing partitions are in
natural bijection with clusters.
\item\label{lat} The noncrossing partitions associated to a Dynkin quiver
$Q$, denoted $\nc_Q$, form a lattice.  
\end{enumerate}

These properties themselves are not our observations.  
We have already mentioned sources for (1).  Statement (2) 
was first established on a 
type-by-type
basis with a computer check for the exceptional types; a proof which does
not rely on the classification of Dynkin diagrams was given by
Brady and Watt \cite{BW2}.  
Our hope was that by setting these properties 
within a new context, we would gain a better understanding of them, and
also of what  
transpires beyond the
Dynkin case.

\medskip

Let $k$ be an algebraically closed field.
Let $Q$ be an arbitrary finite quiver without any oriented cycles.  
Let $\rep Q$ be the category of finite dimensional representations of $Q$.
We refer to exact abelian and extension-closed subcategories of 
$\rep Q$ as {\it wide}.  
The central object of our researches is $\W_Q$, the set of 
finitely generated wide subcategories of 
$\rep Q$. 
There are a number of algebraic objects which are all in bijection one
with another, summarized by the following theorem:

\begin{theorem}  Let $Q$ be a finite acyclic quiver.  Let $\C = \rep Q$.
There are bijections between the following objects.

\begin{enumerate}
\item clusters in the acyclic 
cluster algebra whose initial seed is given by $Q$.
\item isomorphism classes of basic cluster tilting objects in the cluster category 
$\D^b(\C)/(\tau^{-1}[1])$.
\item isomorphism classes of basic exceptional objects in $\C$ which are tilting on their support.
\item finitely generated torsion classes in $\C$.
\item finitely generated wide subcategories in $\C$.
\item finitely generated semistable subcategories in $\C$. 

\end{enumerate}
If $Q$ is Dynkin or extended Dynkin: 
\begin{enumerate}
\item[7.] the noncrossing partitions associated to $Q$.
\end{enumerate}
If $Q$ is Dynkin:
\begin{enumerate}
\item[8.] the elements of the corresponding Cambrian lattice.
\end{enumerate}
\end{theorem}

Some of these results are already known.  
A surjective map from (1) to (2) was constructed in \cite{BMRT}
and a bijection from (2) to (1) in \cite{CK2}, cf. also the
appendix to \cite{BMRT}.  
Those from (2) to (3) to (4)
are well known but we provide proofs, since we could not find a convenient
reference.  
The bijection from (4) to (5) is new.  
The subcategories in (6) are
included among those contained in (5) by a result of \cite{Ki}; 
the reverse inclusion is new.
Bijections from (8) to (1) and from (8) to (7) 
were given in the Dynkin case \cite{Re2}.  Putting these bijections
together yields a bijection from (1) to (7).  A conjectural description
of this bijection was given in \cite{RS}; we prove this conjecture.   
Another bijection between (7) and (8) is also known, though also only
in the Dynkin case
\cite{ABMW}. The extension of the bijection between (1) and (7) to the 
extended Dynkin case is new.

The set $\W_Q$ is naturally ordered by inclusion.   
The inclusion-maximal chains of $\W_Q$
can be identified with the {\it exceptional sequences} for $Q$.  
When $Q$ is of
Dynkin type, $\W_Q$ forms a lattice.  The map from $\W_Q$ to $\NC_Q$ respects
the poset structures on $\W_Q$ and $\NC_Q$, which yields a new proof of the
lattice property of $\NC_Q$ for $Q$ of Dynkin type.    

We also gain some new information about the Cambrian lattices: we 
confirm the conjecture of \cite{Th} that they are {\it trim}, i.e., left modular
\cite{BS}
and extremal \cite{Ma}.  


\begin{acknowledgements} We would like to thank Drew Armstrong, Aslak
Bakke Buan, Fr\'ed\'eric
Chapoton, Matthew Dyer, 
Bernhard Keller, Mark Kleiner, Henning
Krause, Jon McCammond, Nathan Reading, Vic Reiner, Idun Reiten, Claus Ringel, 
Ralf Schiffler, Andrei
Zelevinsky, and an anonymous referee for helpful 
comments and suggestions.
\end{acknowledgements}  

\section{Wide subcategories of hereditary algebras}

\begin{subsection}{Definitions}

In this section we will use some standard facts from homological
algebra, most of which can be found in \cite{bluebook} A.4 and A.5.
In addition to what can be found there we will recall two lemmas.
These facts can be proved with straightforward diagram chases.
The first lemma is a lesser known variant of the snake lemma.

\begin{lemma}\label{wraparound}
If we have maps 
$A \stackrel{\psi}{\rightarrow} B \stackrel{\phi}{\rightarrow} C$ 
in an abelian category then there is a natural exact sequence
\begin{equation*}0 \rightarrow ker \psi \rightarrow ker \phi \psi \rightarrow ker \phi
\rightarrow cok \psi \rightarrow cok \phi \psi \rightarrow cok \phi \rightarrow 0.
\end{equation*}
\end{lemma}

We will also use the fact that pushouts preserve cokernels,
and pullbacks preserve kernels.
\begin{lemma}
Given morphisms $g: A \rightarrow E$ and $f: A \rightarrow B,$
consider the pushout
$$\xymatrix{
A \ar[d]_f \ar[r]^g  & E \ar[d]^{f_*}   \\
B \ar[r]_{g_*}  & E \coprod_A B  
}$$
then $\cok f \simeq \cok f_*$ and $\cok g \simeq \cok g_*$
and the dual statement for pullbacks.
\end{lemma}

Let $k$ be an algebraically closed field.
We will be working with full subcategories of a fixed $k-$linear abelian
category $\C$.
In practice $\C = \rep Q$, the category of finite dimensional modules
over $kQ$ where $Q$ is a finite quiver with no oriented cycles.
In this section we will sometimes prove things in a more general
setting.  We will always assume that $\C$ is small and abelian.
We will also assume that $\C$ has the following three properties:

\begin{description}
\item[Artinian]  Every descending chain of subobjects of an object 
eventually stabilizes.
\item[Krull-Schmidt] Indecomposable objects have local endomorphism rings
and every object decomposes into a finite direct sum of
indecomposables.
\item[Hereditary] The functor $\Ext^1(X,-)$ is right exact for each
object $X$.
\end{description}

The subcategories we consider will always be full and closed under direct sums
and direct summands.  So they are determined by their sets of isomorphism 
classes of indecomposable
objects.  We will abuse notation and occasionally refer to the category
as this set.  Another way of identifying such a subcategory 
is by using a single module.  We let $\add T$
denote the full subcategory, closed under direct sums, 
whose indecomposables are all direct 
summands of $T^i$ for all $i$.  Given a subcategory 
$\A$ of $\C,$ which has only finitely many isomorphism classes of 
indecomposables,
we let $\basic \A$ be the direct sum over a system of representatives
of the isomorphism classes of indecomposables of $\A$.  
So $\add \basic \A =\A$.
We use the operation $\basic$ on a module as shorthand for
$\basic T = \basic \add T$.
Given a full subcategory $\A$ of $\C$ we let $\Gen \A$ be the full 
subcategory whose objects are all quotients of objects of $\A$.
We will also use the same notation $\Gen T$ for an object $T$ in $\C$ as 
shorthand for $\Gen \add T$.

 Some definitions we need for the relevant subcategories include:

\begin{description}
\item[Torsion class] a full subcategory that is
 closed under extensions and quotients.
\item[Torsion free class] a full subcategory that is
 closed under extensions and subobjects.
\item[Exact abelian subcategory]  a full abelian subcategory
where the inclusion functor is exact, hence closed under kernels and cokernels
of the ambient category.
\item[Wide subcategory] an exact abelian subcategory closed 
under extensions.
\end{description}
\end{subsection}

\begin{subsection}{Support tilting modules and torsion classes}
In this section we outline the natural bijection between
basic support tilting modules and finitely generated torsion classes.
We will work in the category $\rep Q$ of finite dimensional
representations
of a finite acyclic quiver $Q$.  
Note that this ambient
category is Artinian, hereditary and satisfies the Krull-Schmidt property.
This material is well known, but we include the results for 
completeness.  Most of the proofs in this section
are given by appropriate references.

\begin{definition}
We say $C$ is a {\it partial tilting module} if 
\begin{enumerate}
\item $Ext^1(C,C) = 0.$
\item $\projd C \leq 1.$
\end{enumerate}
\end{definition}

Note that since we are in a hereditary category the second condition
will always hold.
A {\it tilting module} $C$ is a partial tilting module such that
there is a short exact sequence 
$$ 0 \rightarrow kQ \rightarrow C' \rightarrow C'' \rightarrow 0$$
where $C',C''$ are in $\add C$.

We are particularly concerned with partial tilting modules that
are tilting on their supports.  
For a vertex $x$ in the quiver $Q$,
let $S_x$ be associated simple module of $kQ$.
We say that the support of a module $C$ is the set of simple
modules that occur in the Jordan-Holder series for $C$, up to isomorphism.
This also equals the set of simple modules which occur as subquotients
of finite sums of copies of $C$. 
We need a few lemmas to elucidate the support of a partial 
tilting module.

\begin{lemma}\label{suppstuff}  Let $C$ be a partial tilting module and 
let $M$ be a representation of $Q$.  Then $\supp M \subseteq \supp C$
if and only if $M$ is a subquotient of $C^i$ for some $i$.
\end{lemma}
\begin{proof}  Suppose $\supp M \subseteq \supp C$.  Since the
Jordan-Holder series for $M$ is made up of simples which are
subquotients of $C$, the statement will follow once
we show that 
 the set of subquotients
of $C^i$ for some $i,$ is closed under extension. 
  Suppose that $x,y$ are submodules
of $X,Y$ which are quotients of $C^i$ for some $i$.  We can map an 
extension $e \in \Ext^1(x,y) \rightarrow \Ext^1(x,Y)$, and then since
we are in a hereditary category we can lift via the surjective map
$\Ext^1(X,Y) \rightarrow \Ext^1(x,Y)$ to get an extension $E$ of $Y$ by $X$.
Since $C$ is partial tilting $\Gen C$ is a torsion class closed
under extensions \cite{bluebook}
VI.2.3, so the extension $E$ is in $\Gen C$.  The converse is immediate.
\end{proof}

 A partial tilting module will be called {\it support tilting}
if it also satisfies one of the following equivalent conditions.

\begin{proposition}\label{ptdefns}
The following conditions are equivalent for a partial tilting module $C$.
\begin{enumerate}
\item $C$ is tilting as an $kQ/\ann C$ module.
\item If $M$ is a subquotient of $C^i$ and $Ext^1(C,M)=0$ 
then $M$ is in $\Gen C$.
\item If $\supp M \subseteq \supp C$ and $Ext^1(C,M)=0$
then $M$ is in $\Gen C$.
\item the number of distinct indecomposable direct summands of $C$ is 
the number of distinct simples in its support.
\end{enumerate}
\end{proposition}
\begin{proof}
The equivalence of (1) and (2) is in the proof of
Theorem VI.2.5 in \cite{bluebook}.  The equivalence of 
conditions (1) and (4) follows from Theorem VI.4.4. in \cite{bluebook}.
The equivalence of conditions (2) and 
(3) follows from Lemma~\ref{suppstuff}.
\end{proof}

The following lemma is not used elsewhere, but clarifies the notion
of support tilting.

\begin{lemma}
Suppose that $C$ is a support tilting module.  Then the algebra 
$kQ/\ann C$ is the path algebra of the minimal subquiver
on which $C$ is supported.
\end{lemma}
\begin{proof}
If a vertex $v$ is not in $\supp C$, then clearly the corresponding idempotent
is in $\ann C$ since $e_v C=0$.  Since $\ann C$
is a two sided ideal, any path $x$ that passes through a vertex not
in the support of $C$ is in $\ann C$.  
So this shows that $kQ/\ann C$ is supported on the minimal subquiver $Q'$
on which $C$ is supported.  So we can restrict attention to $Q'$.
Now $C$ is support tilting, and in particular tilting on $Q'$.
Therefore $C$ is faithful by Theorem VI.2.5~\cite{bluebook} and 
so its annihilator is zero on $Q'$.
\end{proof}



We say that an object $P$ in a subcategory $\T$ 
is {\it $\T-$split projective} if all surjective morphisms 
$I \twoheadrightarrow P$ in $\T$ are split.  We say that $P$ is {\it $\T-\Ext$
projective} if $\Ext^1(P,I)=0$ for all $I$ in $\T$.  We will drop the $\T$
in the notation when it is clear from context.    
The proof of the next lemma follows easily from these definitions.

\begin{lemma} \label{extequalssplit}
If the subcategory $\T$ 
is closed under extensions and $U$ is split projective in $\T$, then $U$
is $\Ext$ projective.  
\end{lemma}

We say that a subcategory $\T$ is {\it generated} by  $\P \subseteq \T$  if 
$\T \subseteq \Gen \P$.  We say $\T$ is {\it finitely generated}
if there exists a finite set of indecomposable objects in $\T$
that generate $\T$.  We will use this notion for torsion classes and
wide subcategories.

 We say that $U$ is a
{\it minimal generator}
 if for every direct sum decomposition $U \simeq U' \oplus 
U''$ we have that $U'$ is not generated by $U''$.
We next show that a finitely generated torsion class has a unique minimal
generator.

\begin{lemma} \label{mingen}
A finitely generated torsion class $\T$ has a minimal generator, unique
up to isomorphism, 
which is the direct sum of all its indecomposable split projectives.  
\end{lemma}
\begin{proof}
Since $\T$ is finitely generated, it follows from the Artinian property that
$\T$ has a minimal generator. 
Suppose that $\T$ is finitely generated by the sum of distinct indecomposables
$U=\oplus U_i$ and suppose that $Q$ in $\T$ is an indecomposable split 
projective.
Since $U$ generates, 
we can find a surjection 
$U^i \twoheadrightarrow Q$.  This surjection must split so the
Krull-Schmidt property allows us to conclude that $Q$ is a summand of  $U$.

For the converse, suppose that
 $\T$ is a torsion class with a minimal generator $U$. 
Let $U_0$ be an indecomposable summand of $U$, and consider a
surjection $\rho :E \twoheadrightarrow U_0$ in $\T$.  
We may apply the proof of \cite{bluebook} Lemma IV.6.1 to show that this map
must split.  Therefore $U$ is split projective.
\end{proof}

\begin{lemma}  Let $Q$ be a finite acyclic quiver.  
Let $\T$ be a finitely generated torsion class in 
$\rep Q$ and let $C$ be the direct sum of its
indecomposable $\Ext$-projectives.  Then $C$ is support tilting.
\end{lemma}
\begin{proof}
Let $U$ be the direct sum of the indecomposable split projectives of 
$\T$.  We know by Lemma~\ref{mingen} that $U$ is a minimal
generator of $\T$.  The proof of VI.6.4 \cite{bluebook} shows
that there is an exact sequence 
$$ 0 \rightarrow kQ/\ann U \rightarrow U^i \rightarrow U' \rightarrow 0$$
where $U'$ is $\Ext$-projective in $\T$ and that $U \oplus U'$ is
a tilting module on $kQ/\ann U$.  Then 
Theorem VI.2.5(d) \cite{bluebook} (as noted in the proof of Lemma VI.6.4 
\cite{bluebook}) shows that the $\Ext$-projectives of $\T$
are all summands of $U \oplus U'$.  So $\basic U \oplus U' \simeq \basic C$
and $C$ is support tilting.
\end{proof}

Given a subcategory $\A$ and an object $Q$ of $\C$,
a right $\A$ approximation of $Q$ is a map
$f : B \rightarrow Q$ where $B$ is in $\A$ and any other morphism
from an object in $\A$ to $Q$ factors through $f$.  This is equivalent
to the map $f_* : \Hom(X,B) \rightarrow \Hom(X,Q)$ being surjective
for all $X$ in $\A$.  Basic properties of approximations 
can be found in \cite{AS}.

The next theorem shows that we can recover a basic support tilting object from 
the torsion class that it generates by taking the sum of the indecomposable
$\Ext$-projectives.

\begin{theorem}\label{supptiltextproj}
Let $C$ be a support tilting object.  Then $\Gen C$
is a torsion class and the indecomposable $\Ext$-projectives
of $\Gen C$ are all the indecomposable summands of $C$.  So $\basic C$ is 
the sum of the indecomposable $\Ext$-projectives of $\Gen C.$
\end{theorem}
\begin{proof}
Let $Q$ be an $\Ext$-projective of $\Gen C$.  In particular $Q$ is in $\Gen C$.
Let $f:B \rightarrow Q$ be an $\add C$ right approximation to $Q$.
Since $Q$ is in $\Gen C$ we know that $f$ is surjective.  Apply
the functor $\Hom(C,-)$ to the short exact sequence
$$ 0 \rightarrow \ker f \rightarrow B \rightarrow Q \rightarrow 0$$
to get the exact sequence
$$ \Hom(C,B) \rightarrow \Hom(C,Q) \rightarrow \Ext^1(C,\ker f) 
\rightarrow \Ext^1(C,B).$$
We know $\Ext^1(C,B)=0$ since $C$ is partial tilting and $B$ is in $\add C$.
We also know that the map $\Hom(C,B) \rightarrow \Hom(C,Q)$ is surjective
so $\Ext^1(C,\ker f) =0$.  Also $\ker f$ is a subquotient of 
$C$ so we can conclude that $\ker f \in \Gen C$ since $C$ is support 
tilting.  Now since $Q$ is an $\Ext$-projective in $\Gen C,$ the map $f$
must be split and so $Q$ is in $\add C$.  So any indecomposable $\Ext$-projective
is a direct summand of $C$.  We know that $C$ is $\Ext$-projective
in $\Gen C$ since $\Ext^1(C,C) =0$
and we are in a hereditary category so $C$ can only have $\Ext$-projective 
summands.  
This also shows that $\Gen C$ is a torsion 
class by \cite{bluebook} Corollary VI.6.2.
\end{proof}

\begin{theorem}\label{torssuptilt}
Let $\C = \rep Q$ where $Q$ is a finite acyclic quiver.  Then 
there is a natural bijection between finitely generated 
torsion classes and basic support
tilting objects given by taking the sum of all indecomposable 
$\Ext$-projectives
and its inverse $\Gen$.
\end{theorem}
\begin{proof}
This follows immediately from the above Theorem~\ref{supptiltextproj} 
and Lemma~\ref{mingen}.
\end{proof}
\end{subsection}

\begin{subsection}{Wide subcategories and torsion classes}

We will now define a bijection between finitely generated torsion
classes and finitely generated wide subcategories.  Let $\T$ be a 
torsion class.  The wide subcategory corresponding to it is defined by
taking those objects of $\T$ such that any morphism in $\T$ whose target
is that object, must have its kernel in $\T$.  More explicitly, 
let $\ab(\T)$ be the 
full subcategory whose objects are in the set 
$$\{ B \in \T \,\, | \mbox{ for all }
(g: Y \rightarrow B) \in \T, \ker g \in \T \}$$

\begin{proposition} Let $\T$ be a torsion class.  Then $\ab(\T)$
is a wide subcategory.
\end{proposition}
\begin{proof}  We first show that $\ab(\pc)$ is closed under kernels.
Let $ f: A \rightarrow B$ be a morphism in $\ab(\pc)$.  We know that $\ker f$
is in $\pc$ by the definition of $\ab(\pc)$.  Let 
$i : \ker f \hookrightarrow A$ be the natural injection.
Take a test morphism 
$g : Y \rightarrow \ker f$ in $\pc$.  The composition $i g :
Y \rightarrow \ker f \hookrightarrow A$ is a morphism in $\pc$ with target 
$A$ in $\ab(\pc)$.  So we know that $\ker (i g)$ is in $\pc$,
 but we also know that
$\ker g = \ker (i g)$ since $i$ is injective.  So we can
conclude that $\ker f$ is in $\ab(\pc)$.

Next we show that $\ab(\pc)$ is closed under extensions.  Suppose $A,B$ are 
in $\ab(\pc)$ and let 
$0 \rightarrow A \stackrel{i}{\rightarrow} E 
\stackrel{\pi}{\rightarrow} B \rightarrow 0$ be an 
extension.  Take a test map $g : Y \rightarrow E$ in $\T$.  Using 
Lemma~\ref{wraparound}
for the composition $\pi  g$ we get an induced exact sequence
$$ 0 \rightarrow \ker g \rightarrow \ker(\pi g) 
\stackrel{\psi}{\rightarrow} A.$$
Since $B$ is in $\ab(\pc)$
and $Y$ is in $\pc$ we can conclude that $\ker(\pi g)$ is in $\pc$.
Since $A$ is in $\ab(\pc)$ we can use the map $\psi$ of the above sequence
to conclude that $\ker g$ is in $\pc$.

Lastly we need to show that $\ab(\pc)$ is closed under cokernels.  
We take a morphism $f: A \rightarrow B$ in $\ab(\pc)$.  Write $C$ for
$\cok f$ and let $g: Y \rightarrow C$ be a a test morphism with $Y$ in $\pc$.
Let $\pi : B \rightarrow C$ be the natural surjection.  Note that we know 
that $\ker \pi = \im f$ is in $\pc$ since $\im f$ is a quotient of $A$.
So we form the pullback $Y \prod_C B$, getting an exact sequence
$$ 0 \rightarrow \ker \pi^* \rightarrow Y \prod_C B 
\stackrel{\pi^*}{\rightarrow}  Y \rightarrow 0. $$
Since $\ker \pi^* \simeq \ker \pi$ and $\T$ is closed under extensions,
we see that the pullback $Y \prod_C B$ is in $\T$.
  Now since 
$B$ is in $\ab(\T)$, the map 
$$g^*: Y \prod_C B \rightarrow B$$ has
kernel in $\T$.  So since $\ker g \simeq \ker g^*,$ 
the test map $g$ has kernel in $\T$.
\end{proof}

The map from wide subcategories to torsion classes is described next.
We first need to show that wide subcategories generate torsion classes.

\begin{proposition}
If $\A$ is a wide subcategory of our ambient hereditary category $\C$ , 
then $\Gen \A$ is a torsion class.
\end{proposition}
\begin{proof}
We only need to show that $\Gen \A$ is closed under extensions.
Let $a,b$ be in $\Gen \A$ with surjections $\pi:A \rightarrow a$
and $\rho: B \rightarrow b$ where $A,B$ are in $\A$.
Let 
$$ 0 \rightarrow a \rightarrow e \rightarrow b \rightarrow 0$$
 be an extension.  Since we are in a hereditary category the map
$\pi_* : \Ext^1(b,A) \rightarrow \Ext^1(b,a)$ is surjective.  So we
can choose a lift of the class of the extension above to obtain an
extension 
$$ 0 \rightarrow A \rightarrow E \rightarrow b \rightarrow 0$$
such that the pushout $\pi_*E = E \coprod_A a$ is isomorphic to $e$.
Now we can simply pull back the class of $E$ to an extension 
$\rho^*E=B \prod_b E$ of $B$ by $A$.  Since $\A$ is closed under extensions
we see that $\rho^* E$ 
is in $\A$.  The natural map $\pi_* \rho^* : \rho^*E \rightarrow e$
is surjective since $\cok \rho^* = \cok \rho = 0 = \cok \pi = \cok \pi_*$.
\end{proof}

The next proposition shows that the operations $\ab$ and $\Gen$ are surjective
and injective respectively, and the composition
$\ab \Gen$ gives the identity.  This proposition is more general than we need; 
we will show that once we restrict
to finitely generated subcategories we can obtain a bijection.

\begin{proposition}
If $\A$ is a wide subcategory then $\A=\ab(\Gen \A)$.
\end{proposition}
\begin{proof}
Suppose an object $B$ is in $\A$.  We wish to show that it is in 
$\ab(\Gen \A)$.  So we take a test map $g: y \rightarrow B$ where
$y$ is in  $\Gen \A $.  So there is an surjection $\pi: Y \rightarrow y$
with $Y$ in $\A$.  Then Lemma~\ref{wraparound} shows that there is an
exact sequence 
$$ 0 \rightarrow \ker \pi \rightarrow \ker g\pi \rightarrow \ker g 
\rightarrow 0.$$
Since $g \pi : Y \rightarrow B$ is a map in $\A$ we see that $\ker g\pi$ is 
in $\A$.  So we see that $\ker g$ is in $\Gen \A$ and so $B$ is in 
$\ab(\Gen \A )$.

Now suppose that $b$ is in $\ab(\Gen \A)$.  Since $b$ is in $\Gen \A$,
we can find an surjection $\pi: B \rightarrow b$ with $B$ in $\A$.
Since $b$ is in $\ab(\Gen \A)$ we know that $\ker \pi$ is in $\Gen \A$
and so we can find another surjection $\rho: K \rightarrow \ker \pi$ where
$K$ is in $\A$.  Let $i : \ker \pi \rightarrow B$ be the natural inclusion.
Now we can conclude that $b \simeq \cok  i \rho$
 and $i \rho : K \rightarrow B$
is a map in the wide subcategory $\A$, hence $b$ is in $\A$.
\end{proof}

We need another characterization of the operation $\ab$ in the next proof
so we show we can also define $\ab$ using only kernels of surjective maps
from split projectives of $\T$.

\begin{proposition}\label{surjprop}
Let $\T$ be a finitely generated torsion class in our ambient category $\C$
and define $$\ab_s(\T)=\{ B \in \T \mid \mbox{ for all surjections }
g:(Z \rightarrow B) \in \T \mbox{ with $Z$ split projective, we have } \ker g \in \T \}.$$
Then $\ab(\T) = \ab_s(\T).$

\end{proposition}
\begin{proof}
It is clear that $\ab(\T) \subseteq \ab_s(\T)$ so take $B$ in $\ab_s(\T)$
and a test map $g: Y \rightarrow B$ with $Y$ in $\T$.
We consider the extension 
$$ 0 \rightarrow \ker g \rightarrow Y \rightarrow \im g \rightarrow 0$$
and let $i: \im g \rightarrow B $ be the natural injection.
Since we are in a hereditary category, we know that the induced map
$i^*:\Ext^1(B,\ker g) \rightarrow \Ext^1(\im g,\ker g)$ is surjective
so we can find $Y'$ such that there is a commutative diagram
$$ \begin{CD}
 0 @>>> \ker g @>>> Y @>>> \im g @>>> 0 \\
 @| @| @V{i^*}VV @V{i}VV @| \\          
0 @>>> \ker g @>>> Y' @>>> B @>>> 0
\end{CD}
$$
 with $Y \simeq
\im g \prod_B Y'.$  Now $B$ is in $\T$ so $\cok g \simeq \cok i$ is in $\T$.
So we have a exact sequence
$$ 0 \rightarrow \ker i^* \rightarrow Y \rightarrow Y' \rightarrow 
\cok i^* \rightarrow 0.$$
Now $\ker i^* =\ker i=0 $ and $\cok i^* = \cok i$ is in $\T$ 
and $Y$ is in $\T$ so we may conclude that $Y'$ is in $\T$ since
$\T$ is closed under extensions.  Now we have a surjection 
$g':Y' \rightarrow B$
in $\T$ with kernel isomorphic to $\ker g$.  Let $h:Z\rightarrow Y'$ be a 
surjection, with $Z$ a split projective.  Then $\ker g'h$ is in $\pt$,
by assumption, and by Lemma~\ref{wraparound}, 
$\ker g' \simeq \ker g$ is a quotient of $\ker g'h$,
so it is also in $\pt$.  
Thus, $B$ is in $\ab(\T)$.
\end{proof}

We now are able to prove that we have a bijection from finitely generated
torsion classes to finitely generated wide subcategories.

\begin{proposition}\label{abprogen}
If $\T$ is a finitely generated torsion class then $\ab(\T)$ is 
finitely generated and $\Gen \ab(\T)=\T$.  Furthermore the projectives
of $\ab(\T)$ are the split projectives of $\T$.
\end{proposition}
\begin{proof}
We first show that any $\T$-split projective $U$   
is also in $\ab(\T)$.  Since any surjection $Q \rightarrow U$ in $\T$
splits, and $\T$ is closed under direct summands, we know that $U$ is in
$\ab(\T)$.  Also, since $U$ is $\T$-split projective it is also 
a projective object in $\ab(\T)$.
Conversely, any object $P$ in $\ab(\T)$ admits a surjection from some $U^i$
where $U$ is a split projective generator of $\T$, cf. Lemma~\ref{mingen}.  If
$P$ is projective in $\ab(\T)$, then this surjection must split, so the
projectives of $\ab(\T)$ and the split projectives of $\T$ coincide.  

Now Lemma~\ref{mingen} shows that 
$\T$ is generated by its split projectives,
so we see that $\ab(\T) \subseteq \T$ is also finitely 
generated.
\end{proof}

Combining the above propositions immediately gives one of our main results.

\begin{corollary}  
There is a bijection between finitely generated torsion classes
in $\C$ and finitely generated wide subcategories.  The bijection
is given by $\ab$ and its inverse $\Gen$.
\end{corollary}

\begin{lemma}\label{extprojsubscld}
Let $\C$ be a subcategory of a hereditary category.
If $P$ in $\C$ is $\Ext-$projective, then any subobject $Q \hookrightarrow
P$ in $\C$ is also $\Ext-$projective.  
\end{lemma}
\begin{proof}
If $a$ is in $\C$ then  $\Ext^1(P,a)=0$ and we have a surjection $\Ext^1(P,a)
\twoheadrightarrow \Ext^1(Q,a).$
\end{proof}

\begin{lemma}\label{subsplit}
Let $\T$ be a finitely generated 
torsion class and let $Q$ be a split projective in 
$\T$.  Then any subobject of $Q$ that is in $\T$ is split projective.
\end{lemma}
\begin{proof}
Let $i:P \rightarrow Q$ be an injection 
in $\T$.  Note that $\cok i$ is 
in $\T$.  Since $\T$ is generated by its split projectives we can 
find a surjection $f:R \rightarrow P$ where $R$ is split projective.
Since we are in a hereditary category we can lift the extension $R$
in $\Ext^1(P,\ker f)$ to an extension $E$ in $\Ext^1(Q,\ker f)$.
So we have an exact sequence 
$$ 0 \rightarrow R \rightarrow E \rightarrow \cok i \rightarrow 0$$
which shows that $E$ is in $\T$.  Therefore the surjection $E\rightarrow Q$
must split and the class of $E$ in $\Ext^1(Q,\ker f)$ is zero.  Therefore
the class of $R$ in $\Ext^1(P,\ker f)$ is also zero and so this extension
splits.  So $P$ is a direct summand of the split projective $R$.
\end{proof}

\begin{corollary} If $\pa$ is a finitely generated wide subcategory of
$\rep Q$, then it is hereditary. \end{corollary}

\begin{proof} $\pa=\ab(\Gen \pa)$.  The above result combined with 
Proposition~\ref{abprogen} shows that this category
is hereditary.  \end{proof}

We are also in a position to notice that $\ab(\T) \simeq \rep Q'$
for some finite acyclic quiver $Q'$  
as in the next corollaries.

\begin{corollary}
If $\pa$ is a finitely generated wide subcategory of $\rep Q$,
 then $\pa \simeq \mod \End(U)$
where $U$ is the direct sum of the projectives of $\pa$.
\end{corollary}
\begin{proof} $\pa=\ab(\Gen \pa)$.  Now
Proposition~\ref{abprogen} shows that the abelian category $\pa$ has 
a projective generator which is the sum of the indecomposable 
split projectives in $\Gen \pa$.
So standard Morita theory proves the above equivalence \cite{MR} 3.5.5.  
\end{proof}

\begin{corollary}\label{wideisquiver}
If $\pa$ is a finitely generated wide subcategory of $\rep Q$
then there is a finite acyclic quiver $Q'$ such that
$\ab(\T) \simeq \rep(Q')$.
\end{corollary}

The proof follows on combining the above statements with the theorem 
that a finite dimensional basic hereditary algebra over an algebraically
closed ground field is a path algebra
of an acyclic quiver, \cite{bluebook} Theorem VII.1.7.

We will now proceed to give two alternative characterizations of
the category $\ab(\pt)$.  

\begin{proposition}\label{quot} $\ab(\pt)$ consists of those objects of $\pt$
which can be written as a quotient of a $\pt$-split projective by another
$\pt$-split projective.\end{proposition}

\begin{proof} Suppose $X\in \ab(\pt)$.  Since $\pt$ is generated by split 
projectives, $X$ can be written as a quotient of a split projective.
Now, by the definition of $\ab(\pt)$, the kernel of this map must be in 
$\pt$.  Since it is a subobject of a split projective, it is also a 
split projective.  

Let $X\in \pt$, 
such that $X\simeq P/Q$ for $P,Q$ split projectives.  Let $g:S\rightarrow
X$ be a test morphism, which, 
by Proposition~\ref{surjprop}, we can assume to be surjective, with
$S$ split projective.

From the $\Hom$ long exact sequence, 
we obtain $\Hom(S,P)\rightarrow \Hom(S,P/Q)\rightarrow \Ext^1(S,Q)=0$.  
So $g$ lifts to a map from $S$ to $P$.  We now have a short exact
sequence:
$$0\rightarrow \ker g\rightarrow S\oplus Q \rightarrow P\rightarrow 0$$

Since $P$ is split projective, this splits, and $\ker g$ is a summand of 
$S\oplus Q$, so is in $\pt$.    
So $X\in \ab(\pt)$.  
\end{proof}

We need the following alternative characterization of the category $\ab(\T)$
in the sequel.  It describes $\ab(\T)$ as the perpendicular of the
non-split projectives in $\T$.

\begin{proposition} \label{surjenuf}
Let $\T$ be a finitely generated torsion class and let $P$
be the direct
 sum of a system of representatives of the isomorphism classes of 
indecomposable $\Ext$-projectives which are not split projective.
Then $$\ab(\T) = \{ X \in \T : \Hom(P,X)=0 \} =
\{X \in \T : \Hom(P,X)=\Ext^1(P,X)=0 \}.$$
\end{proposition}

\begin{proof} 
Let $Q$ be a split projective.  
We will begin by showing that there are no non-zero morphisms from
$P$ to $Q$.  Suppose, on the contrary, that $f:P\rightarrow Q$ is 
non-zero.  
Since $\im f$ is a
quotient of $P$, it is in $\pt$, so, since it is a subobject of $Q$,
it is split projective. Thus, the short exact sequence:
$$0\rightarrow \ker f \rightarrow P \rightarrow \im f \rightarrow 0$$
splits, and $P$ has a split projective direct summand, contradicting
the definition of $P$.  

Now suppose we have $X$
in $\ab(\pt)$.  By Proposition~\ref{quot}, $X$ can
be written as $Q/R$, for $Q,R$ split projectives. The $\Hom$ long exact
sequence now gives us:
$$0=\Hom(P,Q)\rightarrow \Hom(P,X)\rightarrow\Ext^1(P,R)=0,$$ so 
$\Hom(P,X)=0$, as desired.

To prove the converse, we need 
to recall briefly the notion of minimal approximations.  
A map $f:R\rightarrow X$ is called {\it right minimal} if any map 
$g:R\rightarrow R$ such that $fg=f$, must be an isomorphism.  A map
that is right minimal and a right approximation (as defined before 
Theorem~\ref{supptiltextproj}) is called a minimal right approximation.

Suppose that $X\in \pt$ and $\Hom(P,X)=0$.  
Let $T$ be the sum of the $\Ext$-projectives of $\pt$.   
Consider the minimal right $\add T$ approximation to $X$; call it 
$k:R\rightarrow X$.  Note that $R$ will not include any non-split projective
summands, since these admit no morphisms to $X$.  Let $K$ be the kernel
of this map.  By the properties of minimal approximation, the map 
$\Hom(T,R)\rightarrow \Hom(T,X)$ is surjective, so $\Ext^1(T,K)=0$.  
Since the support of $K$ is contained in the support of $T$, this implies
that $K$ is in $\pt$ by Proposition~\ref{ptdefns} and 
Theorem~\ref{torssuptilt}.  Since $K$ is a subobject in $\pt$ of a split 
projective, 
$K$ is also split projective.  Now $X\simeq R/K$ shows that $X$ is in 
$\ab(\pt)$, by Proposition~\ref{quot}.  
\end{proof}

A {\it torsion free class} in a category $\mathcal{C}$ 
is the dual notion to a torsion class: it is
a full subcategory closed under direct summands and sums, extensions, and
subobjects.  In the context of representations of a hereditary algebra $A$,
in which, as we have seen, finitely generated wide subcategories are in
bijection with finitely generated torsion classes, it is true dually that
finitely cogenerated wide subcategories are in bijection with finitely 
cogenerated 
torsion free classes.  (Note also that by Corollary~\ref{wideisquiver}
and its dual, finitely cogenerated wide subcategories coincide with 
finitely generated wide subcategories.)
We shall not need to make use of this matter,
so we shall not pursue it here.  

However, we shall need certain facts about torsion and torsion free classes.
These facts are well-known, \cite{bluebook} VI.1.

\begin{lemma}\label{tf} \begin{itemize} \item If $\pt$ is a torsion class in 
$\rep Q$, then the full
subcategory
$\prf$ consisting of all objects admitting no non-zero morphism from an
object of $\pt$, is a torsion free class.  \item Dually, if $\prf$ is a 
torsion free class, then the full subcategory $\pt$ consisting of the 
objects admitting no non-zero morphism to any
object of $\prf$ forms a torsion class.  
\item These operations which construct a torsion free class from a 
torsion class and vice versa are mutually inverse.  Such a pair $(\T,\F)$ of 
reciprocally determining torsion and torsion free classes is called a 
torsion pair.  
\item Given a torsion pair $(\T,\F)$ and an object $X\in \mod A$, there
is a canonical short exact sequence $$0\rightarrow t(X) \rightarrow
X \rightarrow X/t(X) \rightarrow 0$$
with $t(X)\in \pt$ and $X/t(X)\in \prf$.  
\end{itemize}\end{lemma}

\end{subsection}

\begin{subsection}{Support tilting modules and cluster tilting objects}

For $Q$ a quiver with no oriented cycles, the most succinct definition of the
cluster category is that it is $\CL_Q=\pd^b(Q)/\tau^{-1}[1]$, that is to 
say, the bounded derived category of representations of $Q$ modulo 
a certain equivalence.

Fixing a fundamental domain for the action of $\tau^{-1}[1]$, 
we can identify a set of representatives of the isomorphism classes of 
the indecomposable objects of $\CL_Q$ as 
consisting of a copy of the indecomposable representations of $Q$ together 
with 
$n$ objects $P_i[1]$, the shifts of the projective representations.  

A cluster tilting object in $\CL_Q$ is an object $T$ such that 
$\Ext^1_{\CL_Q}(T,T)=0$, and any indecomposable $U$ satisfying
$\Ext^1_{\CL_Q}(T,U)=0=\Ext^1_{\CL_Q}(U,T)=0$ must be a direct summand
of $U$.  
Here $\Ext^j_{\CL_Q}(X,Y)$ is defined as in 
\cite{Ke}, to be $\bigoplus \Ext^j_{\pd^b(Q)}(X,(\tau^{-1}[1])^i(Y))$.

It has been shown \cite{CK2}, cf. also the appendix to \cite{BMRT},
that there is a bijection from
the cluster tilting objects for $\CL_Q$ to the clusters of the acyclic
cluster
algebra with initial seed given by $Q$. The entire structure
of the cluster algebra, and in particular, the exchange relations between
adjacent clusters, can also be read off from the cluster category \cite{BMR}, 
though we shall not have occasion to make use of this here. 

To describe the cluster category $\CL_Q$ 
in a more elementary way, if $X$ and $Y$ are representations of $Q$,
we have that 
$\Ext^1_{\CL_Q}(X,Y)=0$ iff $\Ext^1_{\CL_Q}(Y,X)=0$ iff
$\Ext^1_Q(X,Y)=0=\Ext^1_Q(Y,X)$.   
Additionally, $\Ext^1_{\CL_Q}(X,P_i[1])=0$ iff 
$\Ext^1_{\CL_Q}(P_i[1],X)=0$
iff $\Hom_Q(P_i,X)=0$, and finally, $\Ext^1_{\CL_Q}(P_i[1],P_j[1])=0$ 
always.  
Thus, the condition that an object of $\CL_Q$ is cluster tilting 
can be expressed in terms of conditions that can be checked within $\rep Q$.  

If $T$ is an object in $\CL_Q$, define $\overline T$ to be the maximal direct
summand of $T$ which is an object in $\rep Q$.  
From the above discussion, it is already clear that if $T$ is a cluster
tilting object, then $\overline T$ is a partial tilting object.  
In fact, more is
true: 

\begin{proposition} If $T$ is a cluster tilting object in $\CL_Q$, then
$\overline T$ is support tilting.  Conversely, any support tilting object $V$
can be extended to a cluster tilting object in $\CL_Q$ by adding
shifted projectives in exactly one way.\end{proposition}

\begin{proof}
Let $T$ be a cluster tilting object, which we may suppose to be basic, and
thus to have $n$ direct summands.  
Suppose that $p$ of its indecomposable summands are shifted projectives.
So $\overline T$ has $n-p$ distinct 
indecomposable direct summands.  Observe that the fact that the $p$ 
shifted projective summands have no extensions with $\overline T$ in $\CL_Q$ implies that
$\overline T$ is supported away from the corresponding $p$ vertices of $Q$.  Thus,
$\overline T$ is supported on a quiver with at most $n-p$ vertices.  
But $\overline T$ is a partial 
tilting object with $n-p$ different direct summands, so it must actually
be support tilting.  

Conversely, suppose that $V$ is a support tilting object.  Suppose it has
$n-p$ different direct summands.  Then its support must consist of $n-p$
vertices.  Thus, in $\CL_Q$, the object consisting of the direct sum of 
$V$ and the shifted projectives corresponding to vertices not in the
support of $V$ gives a partial cluster tilting object with $n$ different
direct summands, which is therefore a cluster tilting object.  Clearly, this
is the only way to extend $V$ to a cluster tilting object in $\CL_Q$ by
adding shifted projectives (though there will be other ways to extend $V$
to a cluster tilting object in $\CL_Q$, namely, by adding other 
indecomposable representations
of $Q$).  
\end{proof}

\end{subsection}

\subsection{Mutation}

\newcommand{\cc}{\CL}

An object of $\cc_Q$ is called {\it almost tilting} if it is partial 
tilting and has $n-1$ different direct summands.  A complement to
an almost tilting object $S$ is an indecomposable object $M$ such that
$S\oplus M$ is tilting.  

\begin{lemma}[\cite{BMRRT}] An almost tilting object $S$ in $\cc_Q$ has
exactly two complements (up to isomorphism).  \end{lemma}

The procedure which takes a tilting object and removes one of its 
summands and replaces it by the other complement for the remaining almost
tilting object is called {\it mutation}.  It is the analogue in the
cluster category of the mutation operation in cluster algebras.

Given an object $T$ in $\cc_Q$, we will write $\Gen T$ for the 
subcategory of $\rep Q$ generated by the summands of $T$ which lie in 
$\rep Q$.  When we say that an indecomposable of $T$ is split projective
in $\Gen T$, we imply in particular that it is in $\rep Q$.

The main result of this section is the following proposition:

\begin{proposition}\label{alt} 
If $S$ is an almost tilting object in $\cc_Q$ and
$M$ and $M^*$ are its two complements in $\cc_Q$, then either
$M$ is split projective in $\Gen (M\oplus S)$ or $M^*$ is split
projective in $\Gen(M^*\oplus S)$ and exactly one of these holds.
\end{proposition}

\begin{proof} 
If $S$ contains any shifted projectives, we can remove them and remove the
corresponding vertices from $Q$.  So we may assume that $S$ is almost tilting
in $\rep Q$.  
The main tool used in the proof will be the following fact from
\cite{HU}:

\begin{lemma}[\cite{HU}]\label{hulem}  Let $S$ be an almost tilting object in $\rep Q$.  Then either $S$ is not sincere, in
which case there is only one complement to $S$ in $\rep Q$, or $S$ is sincere,
in which case the two complements to $S$ are related by a short exact 
sequence 

\begin{equation}\label{ss2}0\rightarrow M_1 \rightarrow B \rightarrow M_2\rightarrow 0\end{equation}
with $B$ in $\add S$.\end{lemma}

Suppose first that $S$ is not sincere, and that $M$ is its complement in 
$\rep Q$.  Since $S\oplus M$ is tilting, and therefore sincere, $M$ admits
no surjection from $\add S$.  So $M$ is split projective in 
$\Gen(M\oplus S)$.  On the other hand, the other complement $M^*$ 
to $S$ in $\cc_Q$
is not contained in $\rep Q$, so it is certainly not split projective in
$\Gen(M^*\oplus S)$.  

Now suppose that $S$ is sincere, and that its complements are $M_1$ and 
$M_2$, which are related as in (\ref{ss2}).  Clearly $M_2$ is not 
split projective in $\Gen(M_2\oplus S)$, since it admits a surjection from
$B$.  On the other hand, suppose that there was a surjection $B'\rightarrow
M_1$ with $B'\in \add S$.  The non-zero extension of $M_2$ by $M_1$ would
lift to an extension of $M_2$ by $B'$, but that is impossible since $M_2$
is a complement to $S$.  \end{proof}

An order on basic tilting objects was introduced 
by Riedtmann and Schofield \cite{RieS}. It was later studied 
by Happel and Unger in 
\cite{HU}, in the
context of modules over a not necessarily hereditary algebra. 
Their order is defined in terms of a certain subcategory associated
to a basic tilting object:  
$$\pe(T)=\{M\mid \Ext^i_A(T,M)=0 \textrm{ for }i>0\}.$$
This order on basic tilting objects is defined by 
$S<T$ iff $\pe(S)\subset \pe(T)$.  We recall:

\begin{lemma}
[\cite{bluebook}, Theorem VI.2.5]\label{et}
If $T$ is a tilting object in $\rep Q$, $\pe(T)=\Gen T$.
\end{lemma}

For us, it is natural to consider a partial order on a 
slightly larger ground set, the  
set of tilting objects in $\cc_Q$, 
and to take as our definition that $S\leq T$ iff
$\Gen S\subset \Gen T$.  This is equivalent to considering
the set of all finitely generated torsion classes ordered by inclusion.  
We will show later (in section 4.2) that if $Q$ is a Dynkin quiver, this
order is naturally isomorphic to
the Cambrian lattice defined by Reading \cite{Re}.  

\begin{lemma}\label{ord} Let $T$ be a tilting object in $\cc_Q$, let $X$ be
an indecomposable summand of $T$,
and let $V$ be the tilting object obtained by mutation at $X$.  
If $X$ is split $\Ext$-projective
in $\Gen T$, then $T>V$; otherwise, $T<V$.  
\end{lemma}

\begin{proof} Let $S$ be the almost tilting subobject of $T$ which has
$X$ as its complement, and let $Y$ be the other complement of $S$.
If $X$ is split $\Ext$-projective in $\Gen T$, then, by Proposition
~\ref{alt}, $Y$ is not split $\Ext$-projective in $\Gen V$.  Thus,
$\Gen V$ is generated by $S$, and so $\Gen V\subset \Gen T$.  

On the other hand, if $X$ is not split $\Ext$-projective, then 
$Y$ is, and the same argument shows that $\Gen V \supset \Gen S=\Gen T$.
\end{proof}

In fact, more is true.  It is shown in \cite{HU} that if $T$ and $V$ are
tilting objects in $\rep Q$ related by a single mutation, with, say $T>V$, then
this is a cover relation in the order, that is to say, there is no
other tilting object $R\in \rep Q$ with $T>R>V$.  
The proof in \cite{HU} extends
to the more general setting (tilting objects in $\cc_Q$), but the proof
is not simple and we do not refer again to this result, so we do not 
give a detailed proof here.




\begin{subsection}{Semistable categories}
In this section we show
that any finitely generated wide subcategory of $\rep Q$ is a semistable
category for some stability condition. (A result in the converse direction
also holds, cf.
Theorem~\ref{semistableiswide}.)

Recall that $K_0(kQ)$ is a lattice (i.e., finitely generated 
free abelian group) 
with basis naturally indexed
by the simple modules.  Since the simple modules are in turn indexed
by the vertices we will use the set of vertices
$\{e_i\}$ as a basis of $K_0(kQ)$.  We write $\dimv M$
for the class of $M$ in $K_0(kQ)$.  We know that 
$\dimv M = \sum_i \dim_k M_i e_i$.  The Euler form on $K_0(kQ)$
is defined to be the linear extension of the pairing:
$$\la \dimv M,\dimv N\ra = \dim_k \Hom(M,N)-\dim_k\Ext^1(M,N).$$  
For $\alpha=\sum \alpha_ie_i$ and $\beta= \sum \beta_ie_i$ in $K_0(kQ)$
we have:
$$\la \alpha,\beta\ra =\sum_i \alpha_i \beta_i -\sum_{i\rightarrow
j} \alpha_i \beta_j.$$
The Euler form is generally not symmetric, but we obtain a pairing
on $K_0(kQ)$ by symmetrizing:
$$(\alpha,\beta)=\la\alpha,\beta\ra+\la\beta,\alpha\ra.$$


A stability condition \cite{Ki} is a linear function 
$\theta : K_0(kQ) \rightarrow \Z$. 
A representation $V$ of $Q$ is $\theta$-semistable 
if $\theta ( \dimv(V)) =0$ and if $W \subseteq V$ is 
a subrepresentation then $\theta(\dimv(W)) \leq 0$.  We will abbreviate
$\theta(\dimv(V))=\theta(V)$.
Let $\theta_{ss}$ be the subcategory of representations
that are semistable with respect to $\theta$.

The following theorem is in \cite{Ki}.

\begin{theorem} \label{semistableiswide}
Let $\theta$ be a stability condition.  Then $\theta_{ss}$ is wide.
\end{theorem}

We will need the following easy lemma so we record it here.
\begin{lemma}
Let $\theta$ be a stability condition.  Then $\theta_{ss}$ can
also be described as the representations $V$ such that $\theta(V) =0$
and for all quotients $W$ of $V$, we have that $\theta(W) \geq 0.$
\end{lemma}

Let $T$ be a basic support tilting object with direct summands $T_1,\dots,T_r$.
Since $T$ is support tilting, it is supported on a subquiver $Q'$ of $Q$
with $r$ vertices.  Let us number the vertices on which $T$ is supported
by $n-r+1$ to $n$, and number the other vertices 1 to $n-r$.

Let $d_i$ be the function on $K_0(kQ)$ defined by
$$d_i(\dimv(M))=\la T_i,M \ra =\dim_k\Hom(T_i,M)-\dim_k \Ext^1(T_i,M),$$
for $1\leq i \leq r$.
Let $e_j$ be the function on $K_0(\rep Q)$ defined by $e_j(\dimv(M))=\dim_k
M_j$, that is, $e_j$ is just the $j$th component with respect to the usual
basis.

\begin{theorem} For $T=\oplus_{i=1}^r T_i$ a basic support tilting object, 
the abelian category $\ab(T)=\theta_{ss}$ for $\theta$
satisfying:
$$\theta=\sum_{i=1}^r a_id_i + \sum_{j=1}^{n-r} b_j e_j$$ where
$a_i=0$ if $T_i$ is split projective in $\Gen T$, $a_i>0$ if $T_i$ is
non-split projective, and $b_j<0$.
\end{theorem}

\begin{proof} Suppose $\theta$ is of the form given.  
Let us write $\pt$
for $\Gen T$ and $\pa$ for $\ab(\pt)$.  
First, we will prove some statements about the value of $\theta$ on various
objects in $\pc$, then we will put the pieces together.  

If $X\in \pa$, then $X$ does not
admit any homomorphisms from non-split projectives by 
Proposition~\ref{surjenuf}.  But since $X$ is also in $\pt$, $\Ext^1(T_i,X)=0$
for all $i$.  Thus $\theta(X)=0$. 

If $Y$ is in $\pt\setminus\pa$, then, by Proposition~\ref{surjenuf} again,
$X$ admits some homomorphism from a nonsplit projective.  As before,
$\Ext^1(T_i,X)=0$
for all $i$.  It follows that $\theta(Y)>0$.  

If $Z$ is torsion free, on the other hand, we claim that $\theta(Z)<0$.  
Since $Z$ is torsion free, $\Hom(T_i,Z)=0$ for all $i$.  If $\supp(Z)$
is not contained in $\supp(T)$, then some $b_je_j(Z)<0$, and we are done.  
So suppose that $\supp(Z)\subset \supp(T)$.  We restrict our attention
to the quiver $Q'$ where $T$ is tilting.  Now all we need to do is
show that $\Ext^1(T_i,Z)\ne 0$ for some non-split projective $T_i$.

The torsion free class corresponding to $T$ is cogenerated by $\tau(T)$,
so $Z$ admits a homomorphism to $\tau T_i$ for some $i$. In fact, we can
say somewhat more.  There is a dual notion to split projectives for torsion
free classes, namely split injectives, and a torsion free class is cogenerated
by its split injectives.  So $Z$ admits a morphism to some split injective
$\tau T_i$.  We must show that $T_i$ is a non-split projective.

Now observe that (in $\cc_{Q'}$), $\tau T$ is a tilting object.  Let $S$
be the direct sum of all the $T_j$ other than $T_i$.  So $\tau S$ is almost
tilting.   
By the dual version of Proposition~\ref{alt}, if $V$ is the complement 
to $\tau S$ other than $\tau T_i$ then either $V$ is a shifted projective
or $V$ is non-split injective in $\Cogen \tau S$.  Applying $\tau$, we find
that the complement to $S$ other than $T_i$ is $\tau^{-1} V$.  It follows
that the short exact sequence of Lemma~\ref{hulem} must be

$$ 0\rightarrow \tau^{-1} V \rightarrow B \rightarrow T_i \rightarrow 0$$
where $B$ is in $\add S$. Since $T_i$ admits a non-split surjection from
an element of $\add S$, it must be that $T_i$ is non-split projective. 
The morphism from $Z$ to $\tau T_i$ shows that $\Ext^1(T_i,Z)\ne 0$, so
$\theta(Z)<0$.  

We now put together the pieces.  If $X\in \pa$, then $\theta(X)=0$, while
any quotient $Y$ of $X$ will be in $\pt$, so will have $\theta(Y)\geq 0$.
This implies that $X\in\theta_{ss}$.  Now suppose we have some $V\not\in\pa$.
If $V\in \pt$, $\theta(V)>0$, so $V\not\in\theta_{ss}$.  If $V\not\in \pt$,
$V$ has some torsion free quotient $Z$, and $\theta(Z)<0$, so 
$V\not\in\theta_{ss}$.  Thus $\theta_{ss}=\pa$, as desired.  
%
%
%
\end{proof}

\end{subsection}

\section{Noncrossing partitions}

\subsection{Exceptional sequences and factorizations of the Coxeter element}

For this section, we need to introduce the Coxeter group associated to
$Q$, and the notion of exceptional sequences.  
Let $V = K_0(kQ) \otimes \R$ and recall that $(\alpha,\beta)$ is the
symmetrized Euler form.

A vector $v\in V$ is called a positive root if $( v,v)=2$ and
$v$ is a non-negative
 integral combination of the $e_i$.  To any positive root, there
is an associated reflection $$s_v(w)=w-(v,w) v.$$

Let $W$ be the group of transformations of $V$ generated by these reflections.
$W$ is in fact generated by the reflections $s_i=s_{e_i}$.  The pair 
$(W,\{s_i\})$ 
forms a Coxeter system \cite{H} II.5.1.  

For later use, we recall some facts about {\it reflection functors}.
Let $Q$ be a quiver, and let $v$ be a sink in $Q$.  Let $\widetilde Q$ be
obtained by reversing all the arrows incident with $v$.  Then there is a
functor $R^+_v:\rep Q\rightarrow \rep \widetilde Q$ such that, if we write
$P_v$ for the simple projective module supported at $v$, then $R^+_v(P_v)=0$,
and $R^+_v$ gives an equivalence of categories from the full subcategory
$\Sc$ of $\rep Q$ formed by the objects which do not admit $P_v$ as a direct
summand, to the full subcategory $\widetilde{\Sc}$ of $\rep \widetilde Q$
formed by the objects which do not admit $I_v$ as a direct summand. 
The effect of $R^+_v$ on dimension vectors is closely related to
the simple reflection corresponding to $v$: specifically, if $M$ 
does not contain any copies of $P_v$ as indecomposable summands, then
$\dimv R_v(M)=s_v(\dimv M)$.
Dually, there is a reflection functor $R^-_v$
from $\rep \widetilde Q$ to $\rep Q$.  The functors $R^+_v$ and $R^-_v$
induce mutually inverse equivalences between the full subcategories 
$\Sc$ and $\widetilde{\Sc}$.  The
functor $R^+_v$ is left exact and $R^-_v$ is right exact.

The interaction between reflection functors and torsion pairs
can be described as follows.  

\begin{lemma}\label{tfref} 
Let $Q$ be a quiver with a sink at $v$.  Let $(\T,\F)$ be a torsion pair
where the simple projective $P_v$ is in $\F$.  We apply
the reflection functor $R^+_v$ and write $\tilde{\F}$
and $\widetilde{\T}$ for the images of $\F$ and $\T$ in $\rep \widetilde{Q}$.
Then $\widetilde{\F}$ is a torsion free class
and the indecomposables in its
 corresponding torsion class are the  simple injective $\widetilde{I}_{{v}}$
and the indecomposables of $\widetilde{\T}$.
\end{lemma}
\begin{proof}
Suppose $x$ is in $\widetilde{\F}$ and we have a injection $f:y \rightarrow x$.
If $y$ has $\widetilde{I}_{{v}}$ as a direct summand then so does $x$, but 
$\widetilde{I}_{{v}}$ is not in $\widetilde{\F}$, so this is impossible.
If we apply 
the reflection functor $R^-_v$ we get a morphism
$R^-_v(f):R^-_v (y) \rightarrow R^-_v (x)$.  Let $z$ be its kernel, so we have 
the following sequence exact on the left:
$$0 \rightarrow z \rightarrow R^-_v (y) \rightarrow R^-_v (x)$$

Applying $R^+$, which is left exact, we get:

$$0 \rightarrow R^+_v (z) \rightarrow R^+_v R^-_v(y) \rightarrow 
R^+_vR^-_v(x)$$ 

Noting that since $x$ and $y$ do not have 
$\widetilde I_v$ as a direct summand, $R^+_v R^-_v(f)$ is an injection,
we see that $R^+_v(z)=0$, so $z$ is a sum of copies of $P_v$, and thus
$z\in \F$. 

Now consider the short exact sequence 

$$0 \rightarrow z \rightarrow R^-_v(y) \rightarrow \im (R^-_v(f)) \rightarrow 0$$
Since $\im (R^-_v(f))$ is a subobject of $R^-_v(x)\in \F$, it is also in $\F$.  
Since $\F$ is extension closed, it follows that 
$R^-(y)$ is in $\F$, and thus $y$ is in $\widetilde\F$.
It is clear that $\widetilde{\F}$ is closed under extensions, so it is a 
torsion free
class.  

Now let $x$ be an indecomposable in its associated torsion class.  So 
$\Hom(x,y)=0$ for all $y$ in $\widetilde{\F}$.  Then $\Hom(R^-_v x, R^-_v y)=0$
for all $y$ in $\widetilde\F$ and $\Hom(R^-_v x,P_v)=0.$  Since $P_v$ and 
the indecomposables of $R^-_v \widetilde{\F}$ make up all indecomposables
of $\F$ we see that $R^-_v x$ is in $\T$.  So either $x$ is in 
$\widetilde{\T}$ or $x \simeq \widetilde{I}_{{v}}.$
\end{proof}

A Coxeter element for $W$ is, by definition, the product of the
simple reflections in some order.  We will fix a Coxeter element $\cox(Q)$ to
be the product of the $s_i$ written from left to right 
in an order consistent with the arrows in the quiver $Q$.  (If two
vertices are not adjacent, then the corresponding reflections commute,
so this yields a well-defined element of $W$.)

An object $M\in\rep Q$ is called {\it exceptional} if $\Ext^1(M,M)=0.$
If $M$ is an exceptional indecomposable of $\rep Q$, then $\dimv M$ is a 
positive root.  Thus, there is an associated reflection, $s_{\dimv M}$,
which we also denote $s_M$.  

An exceptional sequence in $\rep Q$ is a sequence $X_1,\dots,X_r$ such
that each $X_i$ is exceptional, and for $i<j$, 
$\Hom(X_j,X_i)=0$ and $\Ext^1(X_j,X_i)=0$.  
The maximum possible length of an exceptional sequence is $n$ since
the $X_i$ are necessarily independent in $K_0(kQ) \simeq \Z^n$. 
An exceptional sequence of length $n$ is called {\it complete}.  
The simple representations of $Q$ taken in any linear order compatible with
the arrows of $Q$ yield an exceptional sequence.  



We recall some facts from \cite{C}.

\begin{lemma}[\cite{C}, Lemma 6] If $(X,Y)$ is an exceptional sequence
in $\rep Q$, there
are unique well-defined representations $R_YX$, $L_XY$ such that
$(Y,R_YX)$, $(L_XY,X)$ are exceptional sequences.\end{lemma}

The objects $R_YX$ and $L_XY$ are discussed in several sources,
for example see \cite{Ru}.  They are called {\it mutations};  
note that mutation has a different
meaning in this context than in the context of clusters.

\begin{lemma}[\cite{C}, p.~124] $$\dimv R_YX= \pm s_Y(\dimv X)$$
$$\dimv L_XY=\pm s_X(\dimv Y)$$ \end{lemma}

\begin{lemma}[\cite{C}, Lemma 8] Let $(X_1,\dots,X_n)$ be a complete 
exceptional sequence.  Then 
$$(X_1,\dots,X_{i-1},X_{i+1},Y,X_{i+2},\dots,X_n)$$ 
is an exceptional sequence iff $Y\simeq R_{X_{i+1}} X_i$.  Similarly,
$(X_1,\dots,X_{i-1},Z,X_i,\dots,X_n)$ is an exceptional sequence iff
$Z\simeq L_{X_i}X_{i+1}.$
\end{lemma}

Let $\pb_n$ be the braid group on $n$ strings, with generators $\sigma_1,
\dots,\sigma_{n-1}$ satisfying the braid relations $\sigma_i\sigma_j=
\sigma_j\sigma_i$ if $|i-j|\geq 2$, and $\sigma_i\sigma_{i+1}\sigma_i=
\sigma_{i+1}\sigma_i\sigma_{i+1}$.  It is straightforward to verify:

\begin{lemma}[\cite{C}, Lemma 9] $\pb_n$ acts on the set of all complete
exceptional
sequences by $$\sigma_i(X_1,\dots,X_n)=(X_1,\dots X_{i-1},X_{i+1},R_{X_{i+1}}
X_i,X_{i+2},\dots, X_n).$$  
\end{lemma} 

We can now state the main theorem of \cite{C}:

\begin{theorem}[\cite{C}, Theorem] The action of $\pb_n$ on complete
exceptional sequences is transitive.\end{theorem}

The next theorem follows from the above results.

\begin{theorem}\label{prodex} If $(E_1,\dots,E_n)$ 
is a complete exceptional sequence
in $\rep Q$, then
$s_{E_1}\dots s_{E_n}=\cox(Q)$.  
\end{theorem}

\begin{proof} By the definition of $\cox(Q)$, the statement is true for the
exceptional sequence consisting of simple modules. 
Now we observe that
the product $s_{E_1}\dots s_{E_n}$ is invariant under the action of the
braid group.  Since the braid group action on exceptional sequences
 is transitive, the theorem
is proved.
\end{proof}

\subsection{Defining noncrossing partitions}

In this section, we introduce the poset of noncrossing partitions.  Let $W$
be a Coxeter group.  Let $T$ be the set of all the reflections of $W$, that is,
the set of all conjugates of the simple reflections of $W$.  

For $w\in W$, define the {\it absolute length} of $w$, written $\ell_T(w)$,
to be the length of the shortest word for $w$ as a product of arbitrary
reflections.  Note that this is not the usual notion of length, which 
would be the length of the shortest word for $w$ as a product of simple
reflections.  That length function, which will appear later, we will denote
$\ell_S(w)$.

Define a partial order on $W$ by taking the transitive closure of the
relations $u< v$ if
$v=ut$ for some $t\in T$ and $\ell_T(v)=\ell_T(u)+1$.  We will use the
notation $\leq$ for the resulting partial order.  This order is called 
{\it absolute order}.  

One can rephrase this definition as saying that $u\leq v$ if there is a
minimal-length expression for $v$ as a product of reflections in which
an expression for $u$ appears as a prefix.  

The noncrossing partitions for $W$ are the interval in this absolute order
between the identity element and a Coxeter element.  (In finite type, the
poset is independent of the choice of Coxeter element, but this is not
necessarily true in general.)  We will write $\NC_Q$ for 
the noncrossing partitions in the Coxeter group corresponding to $Q$ with
respect to the Coxeter element $\cox(Q)$.

Inside $\NC_Q$, for $Q$ of finite type, there is yet another way of
describing the order: for $u,v\in \NC_Q$, we have that
$u \leq v$ iff the reverse
inclusion of fixed spaces holds: $V^v \subseteq V^u$ [BW,Be].

\begin{lemma} $\ell_T(\cox (Q))=n$.
\end{lemma}
\begin{proof} By definition, $\cox(Q)$ can be written as 
a product of $n$ reflections.  We just
have to check that no smaller number will suffice.  
To do this, we use an equivalent definition of $\ell_T$
due to Dyer \cite{Dy}: fix a word for $w$ as a product of simple reflections.  
Then $\ell_T(w)$ is the minimum number of simple reflections you need to
delete from the word to be left with a factorization of $e$.  

It is clear that, if we remove any less than all the reflections from 
$\cox(Q)=s_1\dots s_n$, we do not obtain the identity.  So $\ell_T(\cox(Q))=n$.  
\end{proof}



\begin{lemma} For $\pa$ a finitely generated wide subcategory of $\rep Q$, 
$\cox(\pa)\in \NC_Q$.  \end{lemma}

\begin{proof} The simple objects $(S_1,\dots,S_r)$ in $\pa$ form
 an exceptional sequence in 
$\pa$, so also in $\rep Q$.  Extend it to a complete exceptional sequence in
$\rep Q$.  This exceptional sequence yields a factorization for $\cox(Q)$ as a 
product of $n$ reflections which has $\cox(\pa)$ as a prefix, so 
$\cox(\pa)\in \NC_Q$.  \end{proof}


\begin{lemma}\label{anyex} 
If $(E_1,\dots,E_r)$ is any exceptional sequence for 
$\pa$, then $s_{E_1}\dots s_{E_r}=\cox (\pa)$.  \end{lemma}

\begin{proof} This follows from Theorem~\ref{prodex} applied in $\pa$.  
\end{proof}

\begin{lemma} The map $\cox$ respects the poset structures on $\W_Q$ and 
$\NC_Q$, 
in the sense
that if $\pa\subset \pb$ are finitely generated wide subcategories, then $\cox (\pa)<\cox (\pb)$.  
\end{lemma}

\begin{proof} The exceptional sequence of simples for $\pa$ can be extended
to an exceptional sequence for $\pb$.  
Thus $\cox (\pa)$ is 
a prefix of what is, by Lemma~\ref{anyex}, a minimal-length expression for
$\cox(\pb)$.  So $\cox (\pa)<\cox (\pb)$.\end{proof}

We cannot prove that this map is either injective or surjective in general 
type.  However, in finite or affine type, it is a poset isomorphism, as
we shall proceed to show.  

After this paper was distributed in electronic form, the fact that
$\cox$ is  a poset isomorphism was shown for 
an arbitrary quiver without oriented cycles, 
based on a version of Lemma \ref{gend} below, \cite{IS}.  

\subsection{The map from wide subcategories to noncrossing partitions in 
finite and affine type}

For the duration of this section, we will assume that $Q$ is of finite or
affine type.  

\begin{lemma} Let $\cox (\pa)$ be the Coxeter element for a finite type wide
subcategory of $\rep Q$ of rank $r$.  If $\cox (\pa)$ is written as a product
of $r$ reflections from $T$, then the reflections must all 
correspond to indecomposables of $\pa$. \end{lemma}

\begin{proof} Let $\beta_1,\dots,\beta_r$ be the dimension vectors of the
simple objects of $\pa$.  
Being a finite type Coxeter element, $\cox (\pa)$ has no fixed points in 
the span $\langle \beta_1,\dots,\beta_r\rangle$.  Thus, its fixed subspace
exactly consists of $F_\pa=\cap_i \beta_i^\perp$, and is of codimension $r$.  
A product of $r$ reflections
has fixed space of codimension at most $r$, and if it has codimension 
exactly $r$, then the fixed space must be the intersection of the 
reflecting hyperplanes.  Thus, if $\cox (\pa)=s_{M_1}\dots s_{M_r}$, 
then $\dimv M_j$ must lie in $F_\pa^\perp=\langle 
\beta_1,\dots,\beta_r\rangle$.
The only positive roots in the span $\la \beta_1,\dots,\beta_r\ra$ are the positive roots 
corresponding to indecomposable objects of $\pa$, proving the lemma.
  \end{proof}

Given a subcategory $\A$ of $\C$ we write the perpendicular category as
$${}^\perp\A = \{ M \in \C : \Hom(M,V) = \Ext^1(M,V) =0 \mbox{ for all }
V \in \A \}.$$
If $\A$ is a wide subcategory, so is ${}^\perp A$.  This follows from
Theorem 2.3 of \cite{Sc}, and is easy to check directly.

\begin{theorem} If $Q$ is finite or affine, $\cox$ is an injection.
\end{theorem}

\begin{proof}  Let $\pa$ and $\pb$ be two finitely generated 
wide subcategories of $\rep Q$
such that $\cox (\pa)=\cox (\pb)$.  We may extend an exceptional 
sequence for $\pa$ to one for $\rep Q$, and what we add will be
an exceptional sequence for $\uperp\pa$.   So $\cox(\pa) \cox(\uperp\pa) =
\cox(Q)$. 
Hence it follows that $\cox( {\uperp\pa})=
\cox ({\uperp\pb})$.  Now $\pa$ is of finite or affine type, and it is affine
iff there is an isotropic dimension vector in the span of its dimension
vectors.  Since $\rep Q$ has at most a one-dimensional isotropic subspace,
at most one of $\pa$ or $\uperp\pa$ is of affine type.  Thus without loss of 
generality, we can assume that $\pa$ is of finite type. By assumption,
$\cox \pb=\cox \pa$.  Notice also that $r=\ell_T(\cox (\pa))=
\ell_T(\cox (\pb))$ 
is the rank of $\pb$, so the expression for $\cox (\pb)$ as the product of the
reflections corresponding to the simples of $\pb$ is an expression for
$\cox \pb=\cox \pa$ as a product of $r$ reflections.  By the previous lemma,
the simple objects of $\pb$ must be in $\pa$.  Since the ranks of $\pa$
and $\pb$ are equal, $\pb=\pa$. 
\end{proof}

The argument that $\cox$ is surjective is based on the following lemma:

\begin{lemma} If $Q$ is of finite or affine type and $M_i$ are indecomposable
objects whose dimension vectors are positive roots such that 
$\cox (Q)=s_{M_1}\dots s_{M_n}$, then at least one of the 
$M_i$ is post-projective or pre-injective.  
\end{lemma}

\noindent
Note that any wild type quiver $Q$ with at least three vertices
has tilting objects which are regular (i.e., have no post-projective or
pre-injective summand) \cite{Ringel}.  Since a tilting object yields
an exceptional sequence, and therefore a factorization of $\cox(Q)$, this lemma cannot hold for any such quivers.   

\medskip
\begin{proof}  There is nothing to prove in finite type, since in that case every 
indecomposable is post-projective (and pre-injective).  In affine type,
consider the affine reflection group description of $W$ as a semi-direct
product, $W=W_{fin}\ltimes \Lambda$ 
where $\Lambda$ is a lattice of translations.  The Coxeter
element has a non-zero translation component, since otherwise it would
be of finite order, and we know this is not so because if $M$ is an 
indecomposable non-projective object in $\rep Q$, then 
$\dimv(\tau M)=\cox(Q) \dimv M$ \cite{bluebook} Theorem VII.5.8. 
Since all the regular objects are in finite
$\tau$-orbits, their reflecting hyperplanes are in finite
$\cox (Q)$-orbits.  Thus, they must be parallel to the translation component
of $\cox (Q)$.  Now $\cox (Q)$ 
cannot be written as a product of reflections in hyperplanes
parallel to the translation component of $\cox (Q)$, because such a product would
not have the desired translation component.  
Thus, any factorization of $\cox (Q)$ must
include some factor which is pre-injective or post-projective.    
\end{proof}

\begin{lemma}\label{gend} If $Q$ is of finite or affine type and $\cox (Q)=
s_{M_1}\dots s_{M_n}$, then all the $M_i$ are exceptional.  \end{lemma}

\begin{proof}  There is nothing to prove in the finite type case.  
Fix a specific
$M_i$ which we wish to show is exceptional.  If $M_i$ is post-projective or
pre-injective, we are done.  So assume otherwise.  Then by the previous lemma
there is some $M_j$ with $j\ne i$
which is post-projective or pre-injective.  By braid operations, we may
assume that it is either $M_1$ or $M_n$.  Assume the latter.  Assume
further that $M_n$ is post-projective.  
Now
$\cox(Q)s_M\cox(Q)^{-1}=s_{\tau M}$. 
Conjugating by $\cox(Q)$ clearly
preserves the product, and $\tau$ preserves exceptionality.  Thus, we may
assume that $M_n$ is projective.  Applying reflection functors, we may assume
that $M$ is simple projective.  (In this step, the orientation of $Q$ and
thus the choice of $\cox (Q)$ will change.)  Now let $\pa=\uperp M_n$.  
Note that $\pa$ is isomorphic to the representations of $Q$ with the vertex corresponding
to $M_n$ removed, so $\pa$ is finite type.  Thus, 
$\cox( {\pa})=\cox (Q)s_{M_n}$ is 
a Coxeter element of finite type, so any factorization of it into 
$n-1$ reflections must make use of reflections with dimension vectors in 
$\pa$.  Thus $M_i\in \pa$, so it is exceptional.

If $M_n$ was pre-injective instead of post-projective, we
would have conjugated by $\cox^{-1} (Q)$ to make $M_n$ injective.  The effect
of conjugating by $\cox^{-1} (Q)$ one more time is to turn $s_{M_n}$ into a 
reflection corresponding to an indecomposable projective.  
Then we proceed as above.  
\end{proof}

\begin{theorem} In finite or affine type, the map $\cox$ 
is a surjection.
\end{theorem}

\begin{proof} The argument is by induction on $n$.  Let $w\in \NC_Q$.  
If $w$ is rank $n$, the statement is immediate.  By the previous lemma,
the statement is also true if $w$ is rank $n-1$: we know that $\cox (Q)w^{-1}$ is a reflection
corresponding to an exceptional indecomposable 
object $E$, so $w=\cox( {\uperp E})$.  If 
rank $w<n-1$, there is some $v$ of rank $n-1$ over $w$.  By the above
argument, $v=\cox( {\uperp E})$.  Apply induction to $\uperp E$.  
\end{proof}

\section{Finite type}

Throughout this section, we assume that $Q$ is an orientation of a 
simply laced Dynkin diagram.  A fundamental result is Gabriel's Theorem,
which is proved in \cite{bluebook} VII.5 as well as other sources.
\begin{theorem}
The underlying graph of $Q$ is a Dynkin diagram if and only if 
there is a finite number of isomorphism classes of 
indecomposable representations of $Q$.  In this case $\dimv$ is a bijection
from indecomposable representations of $Q$ to the positive roots in the
root system corresponding to $Q$ expressed with respect to the basis of
simple roots.
\end{theorem}   
In section 4.4, we show how our results extend to 
non-simply laced Dynkin diagrams.  

\subsection{Lattice property of $\NC_Q$}
Our first theorem in finite type is an immediate corollary of results
we have already proved.  This theorem was first established by 
combinatorial arguments in
the classical cases, together with a computer check for the exceptionals.
It was given a type-free proof by Brady and Watt \cite{BW2}.

\begin{theorem} In finite type $\NC_Q$ forms a lattice.
\end{theorem}

\begin{proof} 
If $\pa,\pb\in\W_Q$, then $\pa\cap\pb\in\W_Q$, since 
the intersection
of two abelian and extension-closed subcategories is again
abelian and extension-closed, while the finite generation condition
is trivially satisfied because we are in finite type.  This shows 
that $\W_Q$, ordered by inclusion, has a meet operation.  Since
it also has a maximum element, and it is a finite poset, this suffices
to show that it is a lattice.  Now $\cox$ is a poset isomorphism from
$\W_Q$ ordered by inclusion to $\NC_Q$, so $\NC_Q$ is also a lattice.  
\end{proof}

Note that if $Q$ is not of finite type, $\NC_Q$ need not form a lattice.
(There are non-lattices already in $\widetilde A_n$ for some choices of
(acyclic) orientation \cite{Di}.)  This seems natural
from the point of view of $\W_Q$, since the intersection of two 
finitely-generated subcategories of $\rep Q$ need not be finitely generated.

\subsection{Reading's bijection from noncrossing partitions to clusters}
Our second main finite type 
result concerns bijections between noncrossing partitions
and clusters.  One such bijection in finite type was constructed
by Reading \cite{Re2}, and another subsequently by 
Athanasiadis et al. \cite{ABMW}.  
We will show that the bijection
we have already constructed between clusters and noncrossing partitions
specializes in finite type to the one constructed by Reading.

We first need to introduce Reading's notion of a 
$c$-sortable element of $W$, where $c$ is a Coxeter element for $W$.  
There are several equivalent definitions; we
will give the inductive characterization, as that will prove the most
useful for our purposes.  

A simple reflection $s$ is called {\it initial} in $c$ if there is a reduced
word for $c$ which begins with $s$.  (Note, therefore, that there
may be more than one simple reflection which is initial in $c$, but
there is certainly at least one.)
If $s$ is initial in $c$, then
$scs$ is another Coxeter element for $W$, and $sc$ is a Coxeter element for
a reflection subgroup of $W$, namely, the subgroup generated by the
simple reflections other than $s$.  

By Lemmas 2.4 and 2.5 of \cite{Re2}, and the comment after them, 
the $c$-sortable elements can be characterized by the following properties:
\begin{itemize}
\item  The identity $e$ is $c$-sortable for any $c$.
\item
If $s$ is initial in $c$, then:
\begin{itemize} \item If $\ell_S(sw)>\ell_S(w)$ then $w$ is $c$-sortable iff
$w$ is in the reflection subgroup of $W$ generated by the simple reflections
other than $s$, and $w$ is 
$sc$-sortable.
\item If $\ell_S(sw)<\ell_S(w)$ then $w$ is $c$-sortable iff $sw$ is $scs$-sortable.\end{itemize} \end{itemize}

Let $\Phi$ be the root system associated to $Q$, with $\Phi^+$ the 
positive roots.  
For $w\in W$, we write $I(w)$ for the set of positive roots $\alpha$ such that 
$w^{-1}(\alpha)$ is a negative root.  $I(w)$ is called the {\it inversion set} 
of $w$. 

Gabriel's Theorem tells us that $\dimv$ is a bijection
from indecomposable objects of $\rep Q$ to $\Phi^+$.  
If $\pa$ is an additive subcategory
of $\rep Q$ that is closed under direct summands,
 we write $\Ind(\pa)$ for the corresponding set of positive roots.
If $\alpha\in\Phi^+$, we write $M_\alpha$ for the corresponding 
indecomposable objects.  If $M_\alpha$ is projective (respectively,
injective) we sometimes write $P_\alpha$ (respectively, $I_\alpha$) to
emphasize this fact.  

\begin{theorem} For $Q$ of finite type, there is a bijection between
torsion classes and $\cox (Q)$-sortable elements, $\pt \rightarrow w_\pt$, 
where $w_\pt$ is defined by the property
that $\Ind(\pt)=I(w_\pt)$.
\end{theorem}

\begin{proof} Let $\pt$ be a torsion class.  
We first prove that $\Ind(\pt)$ is the inversion set of some
$\cox (Q)$-sortable element.  
The proof is by induction on the number of vertices
of $Q$ and $|\Ind(\pt)|$.  

Let $\alpha$ be the positive
root corresponding to a simple injective for $Q$.  Let $v_\alpha$
designate the corresponding source of $Q$.  
Now $s_\alpha$ is
initial in $\cox( Q)$. 
If $I_\alpha\not\in \pt$,
then $\pt$ is supported away from $v_\alpha$.  Let $Q'$ be the 
subquiver of $Q$ with $v_\alpha$ removed, and let $W'$ be the corresponding
reflection group.  
Then $\cox( {Q'})=s_\alpha \cox (Q)$ 
and, by induction, $\Ind(\pt)$ is the inversion set of a $\cox({Q'})$-sortable
element $w$.  Now $\ell_S(s_\alpha w)>\ell_S(w)$, and $w$ is 
$s_\alpha \cox(Q)$-sortable, so $w$
is $\cox(Q)$-sortable, as desired.  

Now suppose that $I_\alpha\in\pt$.  In this case, we apply 
the reflection functor 
$R_{v_\alpha}^-$.  Let $\widetilde\pt$ be the image of $\pt$.  
It has one fewer
indecomposable, so by induction, it corresponds to the inversion set of a 
$s_\alpha \cox(Q)s_\alpha$-sortable element, say $\widetilde w$.  
Now $s_\alpha \widetilde w$ 
is $\cox( Q)$-sortable and has 
the desired inversion set.  

Next we show that if $w$ is $\cox(Q)$-sortable then $I(w)$ is $\Ind(\pt)$ for 
some torsion class $\pt$.  
Again, we work by induction.  
If $\ell_S(s_\alpha w)>\ell_S(w)$, then $w$ is $s_\alpha \cox(Q)$-sortable.  
Thus,
by induction, there is a torsion class $\pt'$ on $Q'$ with $\Ind(\pt')=
I(w)$; now $\pt'$ is also a torsion class on $Q$, so we are done.

Suppose on the other hand that $\ell_S(s_\alpha w)<\ell_S(w)$.  
By the induction hypothesis, there is a torsion class
$\widetilde \pt$ on $\widetilde Q$,
with $\Ind(\widetilde\pt)=I(s_\alpha w)$.  Let $\pt$ be the full subcategory
additively generated by $R^+_{v_\alpha}(\widetilde\pt)$ and $I_\alpha$.
Now $\Ind(\pt)=I(w)$.      
By Lemma~\ref{tfref}, $\pt$ is a torsion class.
\end{proof}

The $c$-sortable elements of $W$, ordered by inclusion of 
inversion sets, form a 
lattice, which is isomorphic to the Cambrian lattice $\camb_Q$ \cite{Re3}.  
The reader unfamiliar with Cambrian lattices
may take this as the definition.
(The original
definition of
the Cambrian lattice \cite{Re} involves some lattice-theoretic 
notions which we do not require here, so we shall pass over it.)  Thanks to
the previous theorem, $\camb_Q$ is also isomorphic to the
poset of torsion classes ordered by inclusion.

A {\it cover reflection} of an element $w \in W$ is a reflection 
$t \in T$ such that
$tw=ws$ where $s\in S$ and $\ell_S(ws)<\ell_S(w)$.  

\begin{proposition} If $s$ is initial in $\cox(Q)$, and 
$\pt$ is a torsion class such that $\ell_S(sw_\pt)<\ell_S(w_\pt)$, then
$s$ is a cover reflection for $w_\pt$ 
iff $M_{\alpha_s}$ is in $\ab(\pt)$.  \end{proposition}

\begin{proof} A reflection $t\in T$ corresponding to a positive root
$\alpha_t$ is a cover reflection for $w\in W$ 
iff $I(w)\setminus \alpha_t$ is also the set of inversions for some element of
$W$.  A stronger 
version of the following lemma (without the simply-laced assumption)
is  
\cite{Pil} Proposition 1, see also \cite{Bo} VI\S1 Exercise 16.

\begin{lemma} The sets of roots which arise as inversion sets of
elements of $W$ a simply-laced finite reflection group,
are precisely those whose intersection with any three
positive roots of the form $\{\alpha,\alpha+\beta,\beta\}$ is a subset
which is neither $\{\alpha,\beta\}$ nor $\{\alpha+\beta\}$. 
\end{lemma}

We will say that a set of positive roots is {\it good} if it forms
the inversion set of an element of $W$, and {\it bad} otherwise.  Similarly,
we shall speak of good and bad intersections with a given set 
of positive roots $\{\alpha,\alpha+\beta,\beta\}$.

Thus, if $s$ is not a cover reflection for $w_\pt$, then there are
some positive
roots $R=\{\beta, \beta+\alpha_s,\alpha_s\}$ such that
the intersection of $I(w_\pt)$ with $R$ is good, but becomes bad if we
remove $\alpha_s$.  Thus, $I(w_\pt)\cap R=
\{\beta+\alpha_s,\alpha_s\}$.  So $M_{\beta+\alpha_s}\in\pt$. 
Since $s$ is initial in $c$, we know that $M_{\alpha_s}$ is a simple
injective.  
Thus, there
is a map from $M_{\beta+\alpha_s}$ to $M_{\alpha_s}$, whose kernel will be some
representation of dimension $\beta$.  In fact, though, a generic 
representation of dimension $\beta+\alpha_s$ will be isomorphic to
$M_{\beta+\alpha_s}$ \cite{GR} Theorem 7.1, 
and if we take a generic map from it to 
$M_{\alpha_s}$, 
the kernel will be a generic representation of dimension $\beta$,
thus isomorphic to
$M_\beta$.  Thus,
the kernel of the map from $M_{\beta+\alpha_s}$ to $M_{\alpha_s}$ 
is $M_\beta$. Since $\beta\not\in \Ind(\pt)$, 
we see that $M_\beta\not\in \pt$.  Thus, by the definition of 
$\ab(\pt)$, we have that $M_{\alpha_s}\not\in \ab(\pt)$.  


Conversely, suppose $M_{\alpha_s}\not\in \ab(\pt)$.  
%
%
By Proposition~\ref{surjprop} there is a short exact
sequence  
$0\rightarrow K \rightarrow N \rightarrow M_{\alpha_s}\rightarrow 0$ 
with $K\not
\in \pt$, $N\in\pt$.  Choose such a $K$ so that its total dimension is as small
as possible.  

Let $K'$ be an indecomposable summand of the 
torsion-free quotient of $K$ (as in Lemma~\ref{tf}), 
with respect to the torsion pair $(\pt,\prf)$ determined by $\pt$.    
Then the pushout $N'$ is a quotient of $N$, with 
$0\rightarrow K'\rightarrow N' \rightarrow M_{\alpha_s}\rightarrow 0$.

So by our minimality assumption on $K$, it must be that $K$ is torsion free and
indecomposable.  Suppose $N$ is not indecomposable.  Then let
$N''$ be a direct summand of $N$ which maps in a non-zero fashion to
$M_{\alpha_s}$.  Let $K''$ be the kernel of the map from 
$N''$ to $M_{\alpha_s}$.  Since $K''$ is a subobject of $K$, and $\prf$
is closed under subobjects, by minimality, $K''=K$, so we may assume that 
both $K$ and $N$ are indecomposables, with dimensions, say,
$\beta$ and $\beta+\alpha_s$.  So $\beta\not\in \Ind(\pt)$, 
while $\beta+\alpha_s\in \Ind(\pt)$, as desired.  
%
%
%
%
\end{proof}


Reading's map from $c$-sortable elements to noncrossing partitions can
be characterized by the following proposition:

\begin{proposition}[\cite{Re2}]\label{ncchar} There is a unique map from the $c$-sortable elements
to $\NC_c$ characterized by the properties that $\ncr_c(e)=e$, and,
if $s$ is initial in $c$:
\begin{itemize}
\item
If $\ell_S(sw)> \ell_S(w)$ then $\ncr_c(w)=\ncr_{sc}(w)$.
\item
If $\ell_S(sw)<\ell_S(w)$ and $s$ is a cover reflection of $w$, then $\ncr_c(w)=
\ncr_{scs}(sw)\cdot s$.
\item 
If $\ell_S(sw)<\ell_S(w)$ and $s$ is not a cover reflection of $w$, then
$\ncr_c(w)=s\cdot \ncr_{scs}w\cdot s$ 
\end{itemize}\end{proposition}

There is also a non-inductive definition of the map, but it is 
somewhat complicated, and it will not be needed here, so we do not
give it. 
The above is essentially Lemma 6.5 of \cite{Re2}.


\begin{theorem} The map $\ncr$ coincides with our map from
torsion classes to noncrossing partitions.\end{theorem}

\begin{proof} Our map from torsion classes to noncrossing partitions
is $\cox\circ a$.  The proof amounts to showing that $\cox\circ a$ 
satisfies the 
characterization of Proposition~\ref{ncchar}.  Let $s_\alpha$
be initial in $\cox(Q)$ (and, equivalently, let $M_\alpha$ be a simple 
injective).  
Let $w$ be a $\cox(Q)$-sortable element, and let $\pt$ be the corresponding
torsion class.  
If $\ell_S(s_\alpha w)>\ell_S(w)$, then, as we have seen, $\pt$ is supported on
$Q'$.  The desired condition is now trivially true.  

Now suppose $\ell_S(s_\alpha w)<\ell_S(w)$.  
Define $\widetilde Q$ to be the reflection
of $Q$ at $v_\alpha$.  Let $\widetilde \pt$ be the image of $\pt$ under
the reflection functor $R^-_{v_\alpha}$.  By Lemma~\ref{tfref},  
$\widetilde\pt$ is a torsion class for $\rep \widetilde Q$.  
$\Ind(\widetilde\pt)=s_\alpha(\Ind(\pt)\setminus \alpha)$.  

If $s_\alpha$ is not a cover
reflection for $w$, then $M_\alpha\not\in\ab(\pt)$, so $R^-_{v_\alpha} (\ab(\pt))$ 
is an 
abelian category which generates $\widetilde\pt$, and so $\ab(\widetilde\pt)=
R^-_{v_\alpha}(\ab(\pt))$, and thus $\cox(\ab(\widetilde\pt))=s_\alpha\cox(\ab(\pt))s_\alpha$. 

On the other hand, if $s_\alpha$ 
is a cover reflection for $w$, then $M_\alpha$
is a simple injective for $\ab(\pt)$, and so $R^-_{v_\alpha}$ can be restricted
to a reflection functor for $\ab(\pt)=\rep S$ for some quiver $S$.  Note that
$\ab(\widetilde\pt)$ is contained in $R^-_{v_\alpha} (\ab(\pt))\subset \rep  
\widetilde S$ 
so we can restrict
our attention to the representations of $S$ and $\widetilde S$.  The restriction
of $\pt$ to $\rep S$, though, is all of $\rep S$.  Denote 
the restriction of $\widetilde\pt$ to $\rep \widetilde S$ by $\widetilde \pt _{\widetilde S}$.
Now $\ind \widetilde\pt_{\widetilde S}$ consists of 
all of $\ind \rep \widetilde S$ except $\widetilde M_\alpha$.  
This leaves  us
in a very well-understood situation.  In $\rep \widetilde S$, 
$\widetilde M_\alpha$ is projective,
and if we take $P_{v_\alpha}$ 
to be the projective corresponding to $v_\alpha$ in $\rep S$,
then, in $\rep \widetilde S$, we have that $R^-_{v_\alpha}(P_{v_\alpha})=\tau^{-1}
\widetilde M_\alpha$, so, 
in particular, there is a short exact sequence in  
$\rep \widetilde S$, $0\rightarrow \widetilde M_\alpha
\rightarrow\widetilde P \rightarrow R^-_{v_\alpha}(P_{v_\alpha})\rightarrow 0$, where 
$\widetilde P$ is a 
sum of indecomposable projectives for $\widetilde S$ other than 
$\widetilde M_\alpha$.  
This
shows that $R^-_{v_\alpha}(P_{v_\alpha})$ is not 
split projective.  The other $\Ext$-projectives
of $\widetilde\pt_{\widetilde S}$ are projectives of 
$\rep \widetilde S$, so are
certainly split projectives.  Thus, $\ab(\widetilde\pt_{\widetilde S})$ is the part of
$\rep \widetilde S$ supported away from $\widetilde M_\alpha$, 
and the same is therefore true of
$\ab(\widetilde\pt)$.  Thus, $\cox(\ab(\widetilde\pt))$ can be calculated by taking the
product of the reflections corresponding to the 
injectives of $\rep \widetilde S$ 
other than $s_\alpha$.  
The desired result follows. \end{proof}

Reading also defines a map $\cl_c$ from $c$-sortable elements
to ``$c$-clusters''.  We will
present a version of his map which takes $c$-sortable elements to 
support tilting objects, since that fits our machinery better.  

\begin{proposition}[\cite{Re2}] There is a unique map from $c$-sortable elements to
support tilting objects in $\rep Q$ which can be characterized by the
following properties:
\begin{itemize}
\item If $s$ is initial in $c$ and $\ell_S(sw)>\ell_S(w)$, then 
$\cl_c(w)= \cl_{sc}(w)$.  
\item If $s$ is initial in $c$ and $\ell_S(sw)<\ell_S(w)$, then
$\cl_c(w)= \overline R^+_{v_s}\cl_{scs}(sw)$.
\item $\cl_c(e)=0$.  
\end{itemize}
\end{proposition}

In the above proposition $\overline R^+_{v_s}$ is a map on objects which
is defined by $\overline R^+_{v_s}(T)= R^+_{v_s}(T)$ if $v_s$ 
is in the support of $T$,
but if $v_s$ is not in the support of $T$ then
$\overline R^+_{v_s}(T)=R^+_{v_s}(T)\oplus P_{\alpha_s}$.

\begin{theorem} The map $\cl_c$ corresponds to our map from torsion classes
to support tilting objects.\end{theorem}

\begin{proof} Our map from torsion classes to support tilting objects
consists of taking the $\Ext$-projectives. 
Let $\alpha$ be the positive root corresponding to $s$ initial in $c$, and let
$v$ be the corresponding vertex.  
The image under $R^-_{v}$ of an $\Ext$-projective for 
$\pt$ will be $\Ext$-projective in $\widetilde\pt$.  Conversely, if $M$ is
$\Ext$-projective for $\widetilde\pt$, then $\Ext^1(M,N)=0$ for 
$M,N\in \widetilde\pt$.  It follows that 
$\Ext^1(R^+_{v}(M),R^+_{v}(N))=0$, so in
particular, $\Ext^1(R^+_{v}(M),N')=0$ for $N'$ any indecomposable of
$\pt$ except $M_\alpha$. 
 But $M_\alpha$ is simple injective, so $\Ext^1(R^+_{v}(M),M_\alpha)=0$
as well. The only slight subtlety that can occur is that there might
be an $\Ext$-projective of $\pt$ that is reflected to 0.  (It's not possible
for an $\Ext$-projective of $\widetilde\pt$ to reflect to 0, because $\widetilde\pt$
is by definition the image of $\pt$ under reflection.)  This happens
precisely if $M_\alpha$ is $\Ext$-projective in $\pt$.  

$M_\alpha$ is $\Ext$-projective in $\pt$ iff there are no homomorphisms from
$\pt$ into $\tau(M_\alpha)$, iff there are no morphisms from 
$\widetilde\pt$ into $R^-_{v}(\tau(M_\alpha))$.  
Now $R^-_{v}(\tau(M_\alpha))$ is the injective for
$\rep \widetilde Q$ which corresponds to the vertex $v$.  
There are no morphisms
from $\widetilde\pt$ into $R^-_{v}(\tau(M_s))$ iff $\widetilde \pt $ is
supported away from the vertex $v$.  
\end{proof}

Conjecture 11.3 of \cite{RS} describes the composition $\nc\circ\cl^{-1}$.
An indecomposable $X$ in a support tilting object 
$T$ is {\it upper} if, when we take
$V$ to be the cluster obtained by mutating at $X$, we have that
$\Gen T\supset\Gen V$.  (The definition given in \cite{RS} is not exactly
this, but it is easily seen to be equivalent.)
We can now state and prove the conjecture:

\begin{theorem}[(Conjecture 11.3 of \cite{RS})] For a support tilting object $T$, the fixed space of 
$\cox(\ab(\Gen (T)))$ is the intersection
of the subspaces perpendicular to the roots $\alpha$ corresponding to
upper indecomposables of $T$.\end{theorem}

(Note that, in the finite type setting, it is known that the
fixed subspace of a noncrossing partition determines the noncrossing
partition, so this suffices to describe the map fully.)

\begin{proof}
By Lemma~\ref{ord}, the upper indecomposables of $T$ are exactly
the split $\Ext$-projectives of $\Gen T$.  The fixed space
of $\cox(\ab(\Gen(T)))$ will include the intersection of the 
subspaces perpendicular to the dimension vectors of the split
$\Ext$-projectives, and since 
the fixed
subspace has the same dimension as the intersection of
the perpendicular subspaces, we are done.  
\end{proof}

\subsection{Trimness}

All the lattices which we discuss in this section are assumed to be finite.
An element $x$ of a lattice $L$ is said to be {\it left modular} if, for
any $y<z$ in $L$, 
$$(y\vee x)\wedge z=y\vee(x\wedge z).$$
A lattice is called left modular if it has a maximal chain of
left modular elements.  For more on left modular lattices, see
\cite{BS}, where the concept originated, or \cite{MT}.  

A {\it join-irreducible} of a lattice 
is an element which cannot be written as the join of two strictly smaller
elements, and which is not the minimum element of the lattice.  
A {\it meet-irreducible} is defined dually.  
A lattice is called {\it extremal} if it has the same
number of join-irreducibles and meet-irreducibles as the length of the
longest chain.  (This is the minimum possible number of each.)  See
\cite{Ma} for more on extremal lattices.  

A lattice is called {\it trim} if it is both left modular and extremal.  
Trim lattices have many of the properties of distributive lattices,
but need not be graded.  This concept was introduced and studied in
\cite{Th}, where it was shown 
that the Cambrian lattices in types $A_n$ and $B_n$ 
are trim and conjectured
that all Cambrian lattices are trim.
We will now prove this.  

Let $Q$ be a simply laced Dynkin diagram.  As we have remarked,
the Cambrian lattice $\camb_Q$ can be viewed as the poset of torsion classes
of $\rep Q$ ordered by inclusion, which is the perspective which we shall
adopt.  

The Auslander-Reiten 
quiver for $\rep Q$ is a quiver whose vertices are the isomorphism classes
of indecomposable
representations of $Q,$ and where the number of
 arrows between the vertices associated with indecomposables 
$L$ and $M$ equals the dimension of the space
of irreducible morphisms from $L$ to $M$. 
When $Q$ is Dynkin, this quiver has no oriented cycles.
Thus, one can take a total order on the indecomposables of $Q$ which
is compatible with this order.  We do so, and record our choice by
a map $n:\Phi^+\rightarrow \{1,\dots,|\Phi^+|\}$ so that $n(\alpha)$
records the position of $M_\alpha$ in this total order.  

Let $\ps_i$ be the full additive subcategory, closed under direct summands,
of $\rep Q$ whose indecomposables
are the indecomposables $\{M_\alpha\mid n(\alpha)\geq i\}$.  Each $\ps_i$
is a torsion class.  

\begin{lemma}\label{wedge} For $\pt_1,\pt_2\in\camb_Q$, $\pt_1\wedge\pt_2=\pt_1\cap\pt_2$.
\end{lemma}

\begin{proof}  $\pt_1\cap\pt_2$ is closed under quotients, extensions,
and summands, so it is a torsion class, and thus clearly the maximal
torsion class contained in both $\pt_1$ and $\pt_2$.  \end{proof}


For $\alpha\in\Phi^+$, let $\pt_\alpha=\Gen(M_\alpha)$. 
Recall that  $\Ext^1(M_\alpha,M_\alpha)=0$, so $M_\alpha$ is a partial
tilting object.  Thus, by \cite{bluebook} Lemma VI.2.3,
$\pt_\alpha$ is a torsion
class.  We call such torsion classes {\it principal}.  

\begin{lemma} For $\alpha\in\Phi^+$, the torsion class $\pt_\alpha$
is a join-irreducible in $\camb_Q$.  \end{lemma}

\begin{proof} Let $\pt_\alpha'=\pt_\alpha\cap \ps_{n(\alpha)+1}$.  This
is a torsion class by Lemma~\ref{wedge}, and its indecomposables are
those of $\pt_\alpha$ other than $M_\alpha$ itself.  
Thus, if $\pt_1\vee\pt_2=\pt_\alpha$, then at least one of $\pt_1,\pt_2$
must not be contained in $\pt_\alpha'$, so must contain $M_\alpha$,
and thus all of $\pt_\alpha$.  
\end{proof}

\begin{lemma} The only join-irreducible elements of $\camb_Q$ are the
principal torsion classes.\end{lemma}

\begin{proof} A non-principal torsion class can be written as the
join of the principal torsion classes generated by its split 
$\Ext$-projectives.\end{proof}

\begin{proposition} $\camb_Q$ is extremal.\end{proposition} 

\begin{proof} By the previous lemma,
there are $|\Phi^+|$ join-irreducibles of $\camb_Q$.  Dualizing,
the same is true of the meet-irreducibles.  A maximal chain of
torsion classes $\pt_0\subset\pt_1\subset\dots \subset \pt_m$ must
have $|\pt_{i+1}|\geq |\pt_i|+1$, so the maximal length of such a 
chain is $|\Phi^+|$, proving the proposition.\end{proof}

A torsion class is called splitting if any indecomposable is
either torsion or torsion free.  The $S_i$ are splitting.  

\begin{lemma}\label{vee} If $\ps$ is a splitting torsion class, and $\pt$ is
an arbitrary torsion class, then $\pt\vee\ps=\pt\cup\ps$.  
\end{lemma}

\begin{proof} Let $\prf$ be the torsion free class corresponding to
$\pt$, as in Lemma~\ref{tf}, and let 
$\pe$ be the torsion free class corresponding to $\ps$.  
By the dual of Lemma~\ref{wedge}, $\pe\cap\prf$ is a torsion free
class.  Clearly, the torsion class corresponding to
$\pe\cap\prf$ contains $\ps\cup\pt$.  We claim that equality holds.  
Let $M$ be an indecomposable not contained in $\ps\cup\pt$.
Since $M\not\in \pt$, there 
is an indecomposable $F\in\prf$ which has a non-zero morphism
to $M$.  But since $M\not\in\ps$, $M\in \pe$.  Since $(\ps,\pe)$ forms a 
torsion pair, there are no
morphisms from $\ps$ to $\pe$.  Thus $F$ must not be in $\ps$, and so 
$F\in \pe$, since $(\ps,\pe)$ is splitting.  
We have shown that  $F\in \pe\cap\prf$, and we know there
is a non-zero morphism from $F$ to $M$.  So $M$ is not in the 
torsion class corresponding to $\pe\cap\prf$.  
\end{proof}

\begin{lemma}\label{leftm} Any splitting torsion class is left modular.
\end{lemma}

\begin{proof}
Let $\ps$ be a splitting torsion class.  Let $\pt\supset\pv$ be two
torsion classes.  Now $$\pt\wedge(\ps\vee\pv)=\pt\cap(\ps\cup\pv)
=(\pt\cap\ps)\cup\pv,$$
by Lemmas~\ref{wedge} and~\ref{vee}, and the fact that $\pt\supset \pv$.  
In particular, this implies that $(\pt\cap\ps)\cup\pv$ is a torsion class.
On the other hand,
$\pt\wedge\ps=\pt\cap\ps$.  So $(\pt\wedge\ps)\vee\pv=(\pt\cap\ps)\vee \pv$,
the minimal torsion class containing $\pt\cap\ps$ and $\pv$, which is 
clearly $(\pt\cap\ps)\cup\pv$, as desired.  
\end{proof}

\begin{theorem} $\camb_Q$ is trim.\end{theorem}

\begin{proof} Lemma~\ref{leftm} shows the $S_i$ are left modular, and clearly
they form a maximal chain.  
We have
already showed that $\camb_Q$ is extremal.  Thus, it is trim.  \end{proof}

\subsection{Folding argument}

In our consideration of finite type, we have restricted ourselves to
simply laced cases.  This restriction is not necessary: our conclusions
hold without that assumption.  

The avenue of proof for non-simply laced cases is to apply a {\it folding
argument} in which we consider a simply laced root system which
folds onto the non-simply laced root system.  

Let $Q$ be a simply-laced quiver with a non-trivial automorphism group.  
Define the foldable cluster tilting objects for $Q$ to be those cluster tilting
objects whose isomorphism class is fixed under the action of the
automorphism group of $Q$ on the category of representations,
and similarly for foldable support 
tilting objects.  Define foldable torsion classes of $Q$ to be the
torsion classes of $Q$ stabilized under the action of the automorphism group,
and similarly for foldable wide subcategories.  
Define foldable $c$-sortable elements to be 
those fixed under the action of the 
automorphism group, and similarly for foldable noncrossing partitions. In each
case, the foldable objects for $Q$ 
correspond naturally to the usual object for the 
folded
root system.  All our bijections preserve foldableness, so all our results
go through.  
To conclude that all Cambrian lattices
are trim, we require the fact that the 
sublattice of a trim lattice fixed under a group of lattice automorphisms
is again trim \cite{Th}.

\begin{section}{Example: $A_3$}
In this section we record a few of the correspondences in this paper
for the example of $A_3$ with quiver $Q$: 
$$\xymatrix{
*+[o][F-]{1}  & *+[o][F-]{2} \ar[r]\ar[l]& *+[o][F-]{3}  }$$

The Auslander-Reiten quiver of indecomposable representations
of $Q$ is as follows, where the dimension vectors are written
in the basis given by the simple roots $\alpha_1,\alpha_2,\alpha_3$:

$$\xymatrix{
{[1 0 0]} \ar[dr] & &  
{[0 1 1]} \ar[dr] 
\\
   & {[1 1 1]} \ar[dr]\ar[ur] & &
{[0 1 0]} 
\\
{[0 0 1] \ar[ur]} & &
{[1 1 0] \ar[ur]} 
}$$

In the table on the next page, 
the 14 noncrossing partitions are listed in the
same row as the other objects to which they correspond: the
cluster tilting objects, the support tilting objects, 
the torsion class and the wide
subcategory.  The subcategories of $\rep Q$ are indicated by specifying
a subset of the indecomposables of $\rep Q$, arranged as in the 
Auslander-Reiten quiver.  The support tilting objects and cluster tilting 
objects are indicated
by specifying their summands.    
For the cluster tilting objects, we have drawn a fundamental domain of the 
indecomposable objects in the
cluster category, where the black edges mark the copy of the AR quiver
for $\rep Q$
inside the cluster category, and the dashed edges are maps in 
the cluster category.  The cluster tilting objects can also be viewed
as clusters when the indecomposable objects in the copy of $\rep Q$ are 
identified with positive roots as in the AR quiver above, and the 
three ``extra'' indecomposables are identified with the negative
simple roots $-\alpha_3,-\alpha_2,-\alpha_1$ reading from top to bottom.  
\begin{center}\includegraphics{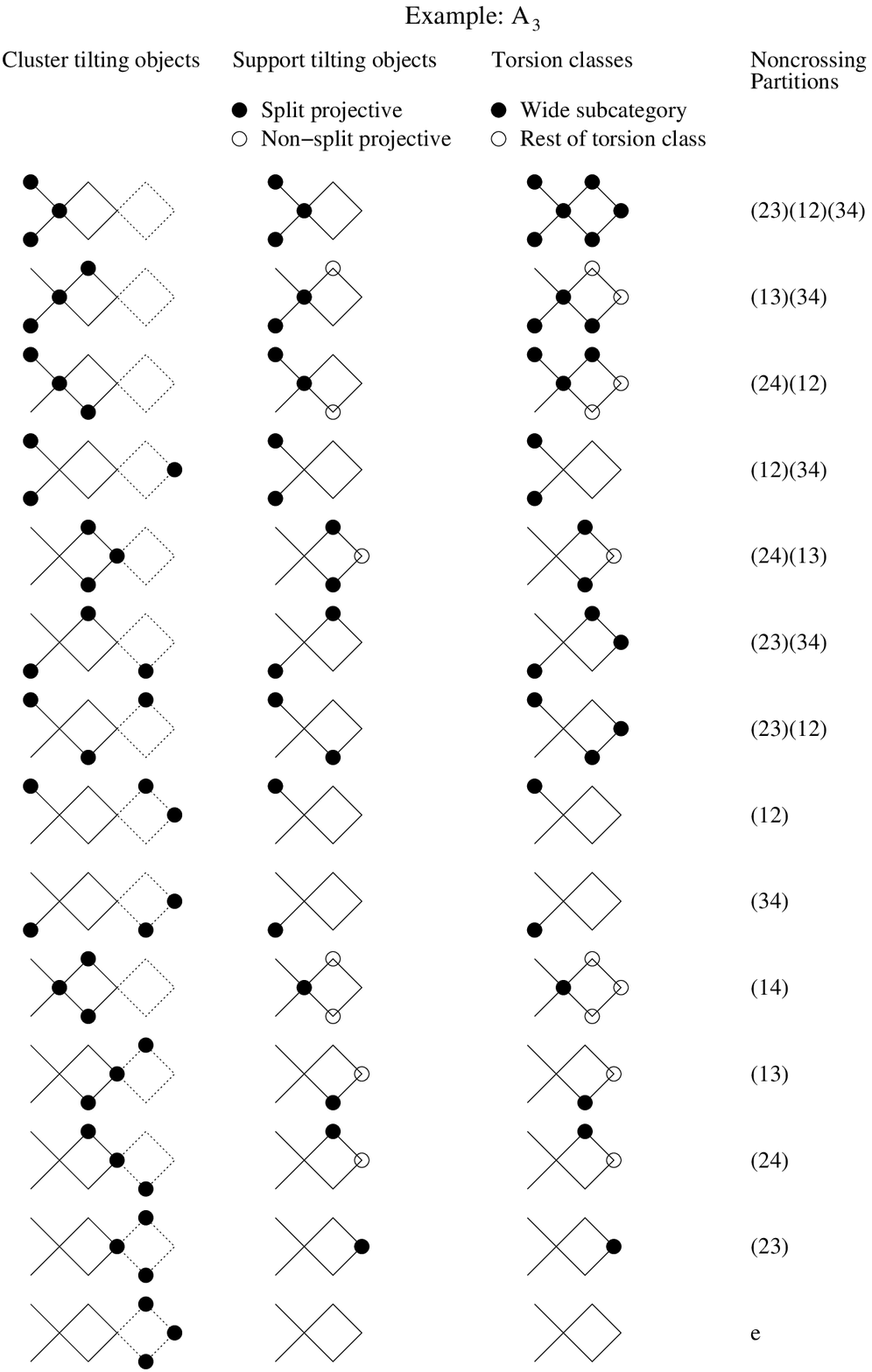}
\end{center}
\end{section}

\end{document}